\magnification=1200

\input amstex

\documentstyle{amsppt}


\hsize=165truemm

\vsize=227truemm


\def\p#1{{{\Bbb P}^{#1}_{k}}}

\def\a#1{{{\Bbb A}^{#1}_{k}}}

\def\Hilb{{{\Cal H}\kern -0.25ex{\italic ilb\/}}}

\def\Scand{{{\Cal S}\kern -0.25ex{\italic cand\/}}}

\def\Hom{{{\Cal H}\kern -0.25ex{\italic om\/}}}

\def\Ext{{{\Cal E}\kern -0.25ex{\italic xt\/}}}

\def\Sim{{{\Cal S}\kern -0.25ex{\italic ym\/}}}

\def\Ker{{{\Cal K}\kern -0.25ex{\italic er\/}}}

\def\PGL{\operatorname{PGL}}

\def\Sing{\operatorname{Sing}}

\def\gr{\operatorname{gr}}

\def\spec{\operatorname{spec}}

\def\proj{\operatorname{proj}}

\def\hom{\operatorname{Hom}}

\def\ext{\operatorname{Ext}}

\def\rk{\operatorname{rk}}

\def\lev{\operatorname{lev}}

\def\Soc{\operatorname{Soc}}

\def\Ofa#1{{{\Cal O}_{#1}}}

\def\M{\operatorname{\frak M}}

\def\mapright#1{\mathbin{\smash{\mathop{\longrightarrow}
\limits^{#1}}}}

\def\mapdown#1{\Big\downarrow\rlap{$\vcenter{\hbox
{$\scriptstyle#1$}}$}}

\def\ga#1{{\accent"12 #1}}


\topmatter

\title
On some Gorenstein loci in $\Hilb_{6}(\p 4)$
\endtitle


\author
Gianfranco Casnati, Roberto Notari
\endauthor
\address
Gianfranco Casnati, Dipartimento di Matematica, Politecnico di Torino,
c.so Duca degli Abruzzi 24, 10129 Torino, Italy
\endaddress

\email
casnati\@calvino.polito.it
\endemail

\address
Roberto Notari, Dipartimento di Matematica, Politecnico di Torino,
c.so Duca degli Abruzzi 24, 10129 Torino, Italy
\endaddress

\email
roberto.notari\@polito.it
\endemail

\keywords
Hilbert scheme, arithmetically Gorenstein subscheme, Artinian algebra
\endkeywords

\subjclassyear{2000}
\subjclass
14C05, 13H10, 14M05
\endsubjclass

\thanks
Both the authors acknowledge support from \lq\lq Programma Vigoni 2005\rq\rq.
\endthanks

\abstract
Let $k$ be an algebraically closed field and let $\Hilb_{d}^{aG}(\p{d-2})$ be the open locus inside the Hilbert scheme $\Hilb_{d}(\p{d-2})$ corresponding to arithmetically Gorenstein subschemes. We prove the irreducibility and characterize the singularities of $\Hilb_{6}^{aG}(\p{4})$. In order to achieve these results we also classify all Artinian, Gorenstein, not necessarily graded, $k$--algebras up to degree $6$. Moreover we describe the loci in $\Hilb_{6}^{aG}(\p{4})$ obtained via some geometric construction. Finally we prove the obstructedness of some families of points in $\Hilb_{d}^{aG}(\p{d-2})$ for each $d\ge6$.
\endabstract

\endtopmatter

\document

\head
1. Introduction and notation
\endhead

Let $k$ be an algebraically closed field and denote by $\Hilb_{p(t)}(\p n)$ the Hilbert scheme parametrizing closed subschemes in $\p n$ with fixed Hilbert polynomial $p(t)\in{\Bbb Q}[t]$. Since A\. Grothendieck proved its existence in [Gr], the problem of its description attracted the interest of many researchers in algebraic geometry and some general and deep properties about such scheme were proved: e.g. we recall its conectedness (see [Ha1]) and a bound on its radius (see [Rv]).

On the other hand, some loci inside the Hilbert scheme were successfully studied, as the locus of codimension $2$ arithmetically Cohen--Macaulay subschemes (see [El]), the locus of codimension $3$ arithmetically Gorenstein subschemes (see [MR] and [K--MR]), the lexicographic point (see [R--S]). In some of the quoted results, it is essential the assumption $\deg(p)\ge1$, because either the case $\deg(p)=0$ needs a different proof or the result does not hold at all. 

However, J\. Fogarty proved that $\Hilb_{d}(\p 2)$ is smooth and irreducible of dimension $2d$ and that, more generally, the same result holds if we consider subschemes of codimension $2$ of a smooth surface (see [Fo]). In [Ia1] the author proved that $\Hilb_{d}(\p n)$ is asymptotically no more irreducible. Indeed for each $d$ and $n$ there always exists a generically smooth component of dimension $nd$ whose general point corresponds to a reduced set of points but, for $d\gg n$, there are also families of larger dimension whose general points correspond to multiple structures of degree $d$ supported on a single point. For a survey on punctual Hilbert schemes see [Ia2].

Since then, there was a big effort in understanding on one hand the structure of the homology ring of $\Hilb_{d}(\p 2)$ (see e.g. [E--S], [I--Y], [Go]), on the other the local properties of some punctual Hilbert schemes $\Hilb_{d}(\p n)$ (see [Ka]). 

From this second viewpoint one checks that a lot of pathologies occur as $n$ increases, so that it is natural to consider, instead of  $\Hilb_{d}(\p n)$, some of its open loci whose closed points parametrize schemes with particular extra technical conditions.

In our paper we are interested in the open locus $\Hilb_{d}^{aG}(\p{d-2})\subseteq\Hilb_{d}(\p{d-2})$ of schemes which are arithmetically Gorenstein. Clearly each point in $\Hilb_{3}^{aG}(\p{1})$ corresponds to a polynomial of degree $3$, thus $\Hilb_{3}^{aG}(\p{1})\cong\p3$. The case $d=4$ has been studied in [A--V] where the more general irreducibility and smoothness of $\Hilb_{4}(\p{2})$ is proved.

Thus it is natural to inspect the case $d\ge5$. When $d=5$, as a by--product, we prove in this paper that $\Hilb_{5}^{aG}(\p{3})$ is irreducible and smooth. We are unable to find other references for this result but we believe that the argument used in  [MR] and [K--MR] for arithmetically Gorenstein subschemes of codimension $3$ and positive dimension in each projective space works also for schemes of dimension $0$ and degree $5$. Such results rest on the well--known beautiful structure Theorem of D\. Buchsbaum and D\. Eisenbud on Gorenstein rings of codimension $3$ (see [B--E]).

For $d=6$ we are able to prove the following

\proclaim{Theorem A}
Assume the characteristic of $k$ is $p\ne2,3$. The locus $\Hilb_{d}^{aG}(\p{d-2})$ is singular for each $d\ge6$. If $d=6$ it is also irreducible and  its singular locus is geometrically described and isomorphic to $\p4$.
\qed
\endproclaim

After [B--E], there was a great attempt of  understanding the structure of Gorenstein rings of codimension at least $4$ and their parameter spaces. In the Artinian graded case it is known, for example, the existence of Gorenstein rings of codimension at least $5$ whose Hilbert function is not unimodal and that in some cases the parameter space is reducible (see [Bo], [G--H--S]).

From this second viewpoint in our paper we classify up to degree $6$ local Artinian $k$--algebras when the characteristic of $k$ is $p\ne2,3$, improving the analogous classification in the locally complete intersection case given in [Ma] and proving the following

\proclaim{Theorem B}
Assume the characteristic of $k$ is $p\ne2,3$. Two different arithmetically Gorenstein, non--degenerate schemes of dimension $0$ and degree $6$ in $\p4$ are projectively equivalent if and only if they are abstractly isomorphic. Furthermore there are exactly $20$ non--isomorphic abstract models for them.
\qed
\endproclaim

Another problem deserving a particular attention is the construction of arithmetically Gorenstein schemes of given dimension and codimension, but it is quite difficult in general to understand if a given scheme can be obtained via a particular construction. In the case of non--degenerate schemes corresponding to closed points in $\Hilb_{6}^{aG}(\p{4})$, we describe the loci obtained via some different constructions.

The structure of the paper is as follows. In Section 2, we observe that the only case when $\Hilb_{d}^{aG}(\p{n})$ contains a point corresponding to a general reduced scheme inside $\Hilb_{d}(\p{n})$ is $n=d-2$ and we recall a result on the uniqueness of the embedding of an arithmetically Gorenstein, non--degenerate scheme of dimension $0$ and degree $d$ in $\p{d-2}$ proved in [C--E]. As a by--product, we compute the minimal free resolution of each such scheme and we prove that no such scheme can be constructed via complete intersection liaison when $d\ge6$.

Section 3 is devoted to the study of the deformations of the most special arithmetically Gorenstein, non--degenerate schemes, which we call {\sl $G$--fat points}\/ in analogy with the usual fat points. A first result we prove is that every other point in $\Hilb_{d}^{aG}(\p{d-2})$ can be flatly deformed to a $G$--fat point. Moreover we show that a $G$--fat point belongs to a component of $\Hilb_{d}(\p{d-2})$ of dimension $d(d-2)$ and its tangent space has dimension $(d^3-7d)/6$, thus proving that $\Hilb_{d}^{aG}(\p{d-2})$ is singular if $d\ge6$ (see Theorem A above).

In Section 4 we classify (non--necessarily graded) Artinian Gorenstein $k$--algebras of degree $d\le6$ for each field $k$ of characteristic $p\ne2,3$, giving the $20$ different abstract models of Theorem B above when $d=6$.

Section 5 contains the proof of the irreducibility and the characterization of the singular locus of $\Hilb_{6}^{aG}(\p{4})$ stated in Theorem A. In Section 6 we describe some different constructions yielding closed points in $\Hilb_{6}^{aG}(\p{4})$ corresponding to non--degenerate schemes. To this purpose it is necessary a case--by--case analysis in order to distinguish which schemes can be obtained via a particular construction.
In the last Section 7, we prove in particular  that the singular locus of $\Hilb_{d}^{aG}(\p{d-2})$ contains many points other the $G$-fat points.

We would like to express our thanks to F. Catanese for some interesting discussions and helpful suggestions. We also thank the referee for her/his careful review.

\subhead
Notation
\endsubhead
In what follows $k$ is an algebraically closed field of characteristic $p\ne2,3$. 

A local ring $R$ is Cohen--Macaulay if
$\dim(R)={\roman{depth}}(R)$. A Cohen--Macaulay ring $R$ is called
Gorenstein if its injective dimension is finite or, equivalently, if $\ext_R^i\big(M,R)=0$ for each $R$--module $M$ and $i>0$.  An
arbitrary ring $R$ is called Cohen--Macaulay (resp. Gorenstein) if
$R_{\frak M}$ is Cohen--Macaulay (resp. Gorenstein) for every
maximal ideal ${\frak M}\subseteq R$. 

All the schemes $X$ are separated and of finite type over $k$. A scheme $X$ is Cohen--Macaulay (resp. Gorenstein) if for each
point $x\in X$ the ring ${\Cal O}_{X,x}$ is
Cohen--Macaulay (resp. Gorenstein). The scheme $X$ is Gorenstein if and only if it is Cohen--Macaulay and its
dualizing sheaf $\omega_{X}$ is invertible.

For each numerical polynomial $p(t)\in{\Bbb Q}[t]$ of degree at most $n$ we denote by $\Hilb_{p(t)}(\p n)$ the Hilbert scheme of closed subschemes of $\p n$ with Hilbert polynomial $p(t)$. With abuse of notation we will denote by the same symbol both a point in $\Hilb_{p(t)}(\p n)$ and the corresponding subscheme of $\p n$.

If $X\in \Hilb_{p(t)}(\p n)$ then the homogeneous ideal $I_X$ of $X$ is saturated and it contains all the homogeneous polynomials in $S:=k[x_0,\dots,x_n]$ vanishing on $X$ and $X=\proj(S_X)$ where $S_X:=S/I_X$ is the projective coordinate ring of $X$. The scheme $X$ is said arithmetically Gorenstein if $S_X$ is a Gorenstein ring.

\head
2. The locus $\Hilb_{d}^{aG}(\p n)$
\endhead

We denote by $\Hilb_{d}^{aG}(\p n)\subseteq\Hilb_{d}(\p n)$ the locus of points representing arithmetically Gorenstein subschemes. For each positive $r$ we set $\Hilb_{d}^{aG,r}(\p n)\subseteq\Hilb_{d}^{aG}(\p n)$ be the subscheme of points $X$ generating a subspace of codimension $r$ in $\p n$. Clearly $\Hilb_{d}^{aG,r}(\p n)=\emptyset$ if $r>n$.

\proclaim{Proposition 2.1}
If $d=1$ then $\Hilb_{d}^{aG}(\p n)=\Hilb_{d}^{aG,n}(\p n)\cong\p n$. If $d=2,3$ then $\Hilb_{d}^{aG}(\p n)=\Hilb_{d}^{aG,n-1}(\p n)$. If $d\ge4$ there is a stratification in disjoint subsets
$$
\Hilb_{d}^{aG}(\p n)=\Hilb_{d}^{aG,n+2-d}(\p n)\cup\Hilb_{d}^{aG,n+1-[d/2]}(\p n)\cup\ldots\cup\Hilb_{d}^{aG,n-1}(\p n).
$$
\endproclaim
\demo{Proof}
The statement in cases $d=1,2,3$ is trivial, thus we assume $d\ge4$ from now on. It is evident that $\Hilb_{d}^{aG,n}(\p n)=\emptyset$.
Let $X\in \Hilb_{d}^{aG,r}(\p n)$ for some $r\le n-1$. With a proper choice of the homogeneous coordinates $x_0,\dots,x_n$ in $\p n$ we can assume that the subspace generated by $X$ is 
$$
\p{n-r}\cong H:=\{\ x_{n-r+1}=\ldots=x_n=0\ \}\subseteq\p n,
$$
thus $x_0,\dots,x_{n-r}$ can be taken as homogeneous coordinates on $H$. If $I\subseteq k[x_0,\dots,x_{n-r}]$ is the homogeneous ideal of $X$ in $H$ then the homogeneous ideal of $X$ in $\p n$ is $I_X:=I+(x_{n-r+1},\dots,x_n)\subseteq k[x_0,\dots,x_n]$. Thus the homogeneous coordinate ring of $X$ in $H$ is $k[x_0,\dots,x_{n-r}]/I\cong k[x_0,\dots,x_n]/I_X= S_X$ which is Gorenstein by assumption, hence we can assume $X\in \Hilb_{d}^{aG}(H)$.

The Hilbert function of such an $X$ then satisfies $h_X(0)=1$ and $h_X(1)=n+1-r$. For each general linear form $\ell\in k[x_0,\dots,x_{n-r}]$, the first difference $\Delta h_X$ is the Hilbert function of the $0$--dimensional graded Gorenstein ring $S_X/\ell$, thus it must be symmetric (see Remark 4.3.9 (a) of [B--V]), and its socle is in some degree $t\ge2$ since $d\ge3$. Moreover $\Delta h_X(0)=1$, $\Delta h_X(1)=n-r$ and its integral is the degree $d$ of the scheme. We conclude that we can have only two different cases: either $\Delta h_X(t)=0$ for $t\ge 3$ whence $n-r+2=d$, or $\Delta h_X(3)\ne0$ whence $n-r+1\le d/2$.
\qed
\enddemo

If $d-2\ge n$ then the description of $\Hilb_{d}^{aG}(\p n)$ can be reduced to the description of the schemes $ \Hilb_{d}^{aG}(\p m)$ with $m<n$, thus the first non--trivial case to examine is $n=d-2$. 
Notice that in this case all the points $\Hilb_{d}^{aG}(\p{d-2})\setminus \Hilb_{d}^{aG,0}(\p{d-2})$ represent schemes contained in a plane, thus complete intersection as soon as $d\le6$. When $d=7,8,9$ this continues to hold true due to Macaulay's growth Theorem (see Theorem 4.2.10 in [B--V]), whereas this is not longer true when $d\ge10$. It is well--known  (see [Sch], Lemma (4.2)) that $d$ general points in $\p{d-2}$ define a point $X\in\Hilb_{d}^{aG,0}(\p{d-2})$, hence $\Hilb_{d}^{aG,0}(\p{d-2})$ is non--empty. 

Since we are interested in dealing with the Hilbert scheme $\Hilb_{d}^{aG,0}(\p{d-2})$, it is helpful to be able to manage flat families in $\Hilb_{d}^{aG,0}(\p{d-2})$. Each such family can be regarded as a flat morphism from a total scheme to a base scheme whose fibres represent points in $\Hilb_{d}^{aG,0}(\p{d-2})$, hence it can be viewed as a Gorenstein cover of degree $d-2$ of a suitable base scheme.

We recall the following result about Gorenstein covers (see [C--E] or [Cs1], Theorem 2.1: in both the papers the Theorem below is wrongly stated without mentioning the condition that the fibres are non--degenerate). 

\proclaim{Theorem 2.2}
Let $X$ and $Y$ be schemes, $Y$ integral and
let
$\varrho\colon X\to Y$ be  a Gorenstein cover of degree $d\ge3$.
There exists a unique $\p{d-2}$--bundle $\pi\colon {\Bbb
P}\to Y$ and an embedding
$i\colon X\hookrightarrow{\Bbb P}$ such that
$\varrho=\pi\circ i$ and $X_y:=\varrho^{-1}(y)\subseteq{\Bbb
P}_y:=\pi^{-1}(y)\cong\p{d-2}$ is a non--degenerate arithmetically Gorenstein subscheme for each
$y\in Y$. Moreover the following properties hold.

\item{i)} ${\Bbb P}\cong{\Bbb P}({\Cal E})$ where ${\Cal
E}$ is the kernel of the trace map $\varrho_*\omega_{X\vert Y}\to\Ofa Y$ (which is dual of the natural morphism $\Ofa Y\to\varrho_*\Ofa X$).

\item{ii)} The composition
$\varphi\colon\varrho^*{\Cal E}\to\varrho^*\varrho_*\omega_{X\vert Y}\to\omega_{X\vert Y}$
is surjective and the ramification divisor $R$ satisfies $\Ofa
X(R)\cong\omega_{X\vert Y}\cong\Ofa X(1):=i^*\Ofa{{\Bbb P}({\Cal E})}(1)$.

\item{iii)} There exists an exact sequence ${\Cal N}_*$ of locally
free
$\Ofa{\Bbb P}$--sheaves
$$ 0\longrightarrow{\Cal N}_{d-2}(-d)\mapright{\alpha_{d-2}}{\Cal
N}_{d-3}(-d+2)\mapright{\alpha_{d-3}}\dots\mapright{\alpha_2}{\Cal
N}_1(-2)\mapright{\alpha_1}\Ofa{\Bbb P}\longrightarrow\Ofa X
\longrightarrow0
\tag 2.1.1
$$
unique up to unique isomorphisms and whose restriction to the
fibre
${\Bbb P}_y$ over $y$ is a minimal free resolution of
the structure sheaf of $X_y$, in particular ${\Cal
N}_i$ is fibrewise trivial. ${\Cal N}_{d-2}$ is invertible and for
$i=1,\dots,d-3$ one has
$$
\rk{\Cal N}_i=\beta_i={{i(d-2-i)}\over{d-1}}{d\choose{i+1}}.
$$
Moreover $\pi^*\pi_*{\Cal N}_*\cong{\Cal N}_*$ and
$\Hom_{{\Bbb P}}\big({\Cal N}_{*},{\Cal N}_{d-2}(-d)\big)\cong{\Cal
N}_{*}$.

\item{iv)} If ${\Bbb P}\cong{\Bbb P}({\Cal E}')$ then ${\Cal
E}'\cong\Cal E$ if and only if ${\Cal N}_{d-2}\cong\pi^*\det{\Cal
E}'$ in the resolution (2.1.1) computed with respect to the
polarization $\Ofa{{\Bbb P}({\Cal E}')}(1)$.
\endproclaim

The above Theorem yields the following two corollaries

\proclaim{Corollary 2.3}
Let $X$ be a Gorenstein scheme of dimension $0$ and degree $d$. Then there exists a non--degenerate embedding $i\colon X\hookrightarrow\p{d-2}$ as arithmetically Gorenstein subscheme. In particular we have a resolution
$$
0\longrightarrow S(-d)\longrightarrow
S(-d+2)^{\oplus\beta_{d-3}}\longrightarrow\dots\longrightarrow
S(-2)^{\oplus\beta_1}\longrightarrow S\longrightarrow S_X\longrightarrow0
\tag 2.3.1
$$
where $S:=k[x_0,\dots,x_{d-2}]$ and
$$
\beta_h:={{h(d-2-h)}\over{d-1}}{d\choose{h+1}},\qquad h=1,\dots,d-3.
$$

Moreover if $j\colon X\hookrightarrow\p{d-2}$ is another embedding whose image is non--degenerate and arithmetically Gorenstein then there exists $\varphi\in\PGL_{d-1}$ such that $\varphi\circ i=j$.
\endproclaim
\demo{Proof}
The first part follows from Theorem 2.2, taking $Y=\spec(k)$: in this case $i$ embeds $X$ inside $\p{d-2}={\Bbb P}({\Cal E})$, ${\Cal E}\cong k^{\oplus d-1}$ being the kernel of the trace morphism $\varrho_*\omega_{X\vert \spec(k)}\to k$.

For the second part we recall that each map $i$, embedding $X$ in $\p{d-2}$ as arithmetically Gorenstein subscheme, is induced by the surjection $\varrho^*{\Cal E}\to\varrho^*\varrho_*\omega_{X\vert \spec(k)}\to\omega_{X\vert \spec(k)}$, thus such an embedding is unique up to automorphisms of $\p{d-2}={\Bbb P}({\Cal E})$.
\qed
\enddemo

\proclaim{Corollary 2.4}
Let $B$ an integral scheme and ${\Cal X}\subseteq\a n\times B$ a closed subscheme such that the induced morphism ${\Cal X}\mapright\varrho B$ has Gorenstein fibres of dimension $0$ and degree $d$.

Then for each $b_0\in B$ there exists an open neighborhood $B'$ of $b_0$ and a morphism $B'\to\Hilb_{d}(\p{d-2})$ whose image is contained in $\Hilb_{d}^{aG,0}(\p{d-2})$ and such that ${\Cal X}':=\varrho^{-1}(B')$ is the pull--back of the universal family over $\Hilb_{d}(\p{d-2})$.
\endproclaim
\demo{Proof}
Since the fibres of $\varrho$ have dimension $0$, the embedding ${\Cal X}\subseteq\a n\times B\subseteq\p n\times B$ implies that $\varrho$ is flat (see [Ha2], Theorem III.9.9) and projective. Since it is quasi--finite then it is actually finite (see [Ha2], Exercise III.11.2). In particular we can interpret $\varrho$ as a Gorenstein cover of degree $d$ of the scheme $B$. Thus we have an embedding $i\colon{\Cal X}\hookrightarrow {\Bbb P}$ in a $\p{d-2}$--bundle such that $\varrho=\pi\circ i$ and the fibres of $\varrho$ are arithmetically Gorenstein subschemes of the fibres of $\pi$. 

Let $B'\subseteq B$ be an open neighborhood of $b_0$ over which $\Bbb P$ is isomorphic to $B'\times \p{d-2}$. Then the universal property of the Hilbert scheme yields the existence of a $B'\to \Hilb_{d}(\p{d-2})$ whose image is in $\Hilb_{d}^{aG,0}(\p{d-2})$.
\qed
\enddemo

By the above results it follows our interest about the locus  $\Hilb_{d}^{aG}(\p{d-2})$. We begin its inspection with the following

\proclaim{Proposition 2.5}
$\Hilb_{d}^{aG,r}(\p{d-2})\subseteq \Hilb_{d}(\p{d-2})$ is open and non--empty if $r=0$, it is constructible if $r\ge1$. $\Hilb_{d}^{aG}(\p{d-2})$ is constructible and it contains a closed irreducible component $\Hilb_{d}^{gen}(\p{d-2})$ of dimension $d(d-2)$
\endproclaim
\demo{Proof}
First of all we notice that $\Hilb_{d}^{aG,0}(\p{d-2})$ is non--empty since it contains the schemes formed by $d$ points in general position in $\p{d-2}$ (see [Sch], Lemma (4.2)).

For each $X\in\Hilb_{d}(\p{d-2})$ its Hilbert function $h_X$ is bounded from above by the function $h_{max}$ defined by $h_{max}(0)=1$, $h_{max}(1)=d-1$, $h_{max}(t)=d$ for $t\ge2$. Due to the lower semicontinuity of the Hilbert function proved in Proposition 1.5 of [B--G] we obtain the existence of a non--empty subscheme ${\Cal U} \subseteq\Hilb_{d}(\p{d-2})$ whose points correspond to schemes $X$ with Hilbert function $h_X=h_{max}$. Sequence (2.3.1) yields that $\Hilb_{d}^{aG,0}(\p{d-2})\subseteq\Cal U$.

The Betti numbers of the schemes in $\Cal U$ are upper semicontinuous by Proposition 2.15 of [B--G]. Since for each such $X$ there is a minimal free resolution of the homogeneous coordinate ring $S_X:=S/I_X$, $S:=k[x_0,\dots,x_{d-2}]$, of the form
$$
0\to F_{d-2}\to\dots\to F_{1}\to S\to S_X\to 0
$$
with $F_{d-2}\ne0$, it follows that there exists an open non--empty subset ${\Cal V}\subseteq\Cal U$ whose points $X$ satisfy $\rk(F_{d-2})=1$, i.e. they are arithmetically Gorenstein. Thus we infer that  $\Hilb_{d}^{aG,0}(\p{d-2})={\Cal V}$ from Corollary 2.3, whence follows the openness of $\Hilb_{d}^{aG,0}(\p{d-2})$ inside $\Hilb_{d}(\p{d-2})$. In a similar way one proves that $\Hilb_{d}^{aG,r}(\p{d-2})$ is locally closed inside $\Hilb_{d}(\p{d-2})$ when $r\ge1$.

Let $\Sim^d(\p{d-2})$ be the $d$--fold symmetric product of $\p{d-2}$ and let ${\Cal W}\subseteq\Sim^d(\p{d-2})$ be the $d(d-2)$--dimensional irreducible open subset parametrizing cycles of $d$ distinct points. For each $d\ge3$ the natural morphism $s\colon \Hilb_{d}(\p{d-2})\to\Sim^d(\p{d-2})$ (see [Fo]) induces an isomorphism $s^{-1}({\Cal W})\to{\Cal W}$. Then $\Hilb_{d}^{gen}(\p{d-2})=\Hilb_{d}^{aG}(\p{d-2})\cap \overline{s^{-1}({\Cal W})}$ is non--empty. 
\qed
\enddemo

\remark{Remark 2.6}
The locus $\Hilb_{d}^{aG}(\p{n})$ does not contain a point corresponding to a general reduced scheme $\Hilb_{d}(\p{n})$ if $n\ne d-2$. The Hilbert function of a general reduced scheme $X$ of degree $d$ in $\p n$ is
$$
h_X(t)=\min\left\{ {n+t\choose n},d\right\}
$$
because general points impose independent conditions on polynomial of any degree. Since $X$ is arithmetically Gorenstein then its first difference $\Delta h_X(t)$ must be symmetric (see the proof of Proposition 2.1). If $n>d-2$ then $\Delta h_X(0)=1$, $\Delta h_X(1)=d-1$, $\Delta h_X(t)=0$ when $t\ge2$, hence $X$ is not arithmetically Gorenstein. Let $n\le d-2$: if $t_0\ge1$ is the greatest integer such that $h_X(t)={n+t\choose n}< d$, then
$$
\gather
1=\Delta h_X(0)=\Delta h_X(t_0+1)=d-{n+t_0\choose n},\\
n=\Delta h_X(1)=\Delta h_X(t_0)={n+t_0\choose n}-{n+t_0-1\choose n}={n\over {t_0}}{n+t_0-1\choose n}.
\endgather
$$
Thus $n=d-1-t_0$. By substituting in the first equality above we finally obtain $n=d-2$. 
\endremark
\medbreak

We conclude the section with the following interesting remark

\remark{Remark 2.7}
One could hope that each $X\in \Hilb_{d}^{aG,0}(\p{d-2})$ is in the same liaison class of a point as in the cases $d=3,4,5$. This does not hold true when $d\ge6$, also in the simplest case $X=\{\ P_0,\dots,P_{d-1}\ \}\in\Hilb_{d}^{gen}(\p{d-2})$, $P_0,\dots,P_{d-1}$ being general points in $\p{d-2}$. Indeed if such an $X$ would linked to a single point,  since the canonical module of $S_X$ is $S_X(1)$ due to Sequence (2.3.1), then $H^0_{(x_0,\dots,x_{d-2})}\big(I_X/I_X^2\otimes S_X(1)\big)=0$ (see [K--M--MR--N--P], Corollary 6.6 and the reference cited there). Let $\Im_X$ be the sheaf of ideals of $X$ in $\Ofa X$. It follows from the exact sequence of graded module over $S$
$$
\align
0&\longrightarrow H^0_{(x_0,\dots,x_{d-2})}\big(I_X/I_X^2\otimes S_X\big)\longrightarrow I_X/I_X^2\mapright\mu \\
&\phantom{\longrightarrow}\longrightarrow \bigoplus_{t\in\Bbb Z} H^0\big(X,\Im_X/\Im_X^2(t)\big)\longrightarrow H^1_{(x_0,\dots,x_{d-2})}\big(I_X/I_X^2\otimes S_X\big)\longrightarrow0
\endalign
$$
that such a vanishing of the local cohomology is equivalent to the injectivity of $\mu$ in each degree. Since $X$ is the union of $d$ simple points in $\p{d-2}$ then $h^0\big(X,\Im_X/\Im_X^2(3)\big)=d(d-2)$. On the other hand taking the degree three component of the exact sequence 
$$
0\longrightarrow I_X/I_X^2\longrightarrow S/I_X^2\longrightarrow S_X\longrightarrow 0,
$$
since $I_X^2$ is generated in degree $4$, we finally obtain that $(I_X/I_X^2)_{3}$ has dimension $(d^3-7d)/6$ which exceed $d(d-2)$ as soon as $d\ge6$.

Since the property of being linked via a complete intersection is stable by generalization, then the above computations also show that no schemes $X\in \Hilb_{d}^{aG,0}(\p{d-2})$ in the same connected component of a reduced scheme can be linked to a complete intersection.

An analogous remark holds for every arithmetically Gorenstein subscheme of degree  $d$ and codimension $d-2$ in $\p n$. Indeed one can reduce to the case described above by taking a general subspace of dimension $d-2$.
\endremark
\medbreak

\head
3. The $G$--fat point and its deformations
\endhead

Let $X=\proj(S_X)\in\Hilb_{d}^{aG,0}(\p{d-2})$. It is non--degenerate by definition and its Hilbert function is $h_X(0)=1$, $h_X(1)=d-1$, $h_X(t)=d$ for $t\ge2$. There exists a hyperplane $H\subseteq\p{d-2}$ such that $X\cap H=\emptyset$. Assuming $H=\{\ x_0=0\ \}$ then $x_0$ is a regular element in $S_X$, whence $R:=S_X/(x_0)$ turns out to be an Artinian Gorenstein graded local $k$--algebra of degree $d$ with maximal ideal ${\frak M}:=(x_1,\dots,x_{d-2})$, thus its components $R_i$ in degree $i$ satisfy $\dim_k(R_0)=1$, $\dim_k(R_1)=d-2$, hence $\dim_k(R_2)=1$ and $R_i=0$ for $i\ge3$.

Since $R$ is Gorenstein then a proper choice of the generators of $\frak M$ allows us to assume that the matrix of the bilinear form 
$$
\psi_1\colon \M/\M^{2}\times\M/\M^{2}\to\M^2\subseteq\Soc(R)
$$
induced by the multiplication is the identity up to scalars,  due to the following well--known

\proclaim{Proposition 3.1}
Let $A:=\oplus_{h=0}^eA_h$ be a graded, local, Artinian $k$--algebra with maximal ideal ${\frak M}:=\oplus_{h=1}^eA_h$, $A_0=k$ and $A_e\ne0$. Then $A$ is Gorenstein
if and only if the forms
$$
\psi_i\colon \M^i/\M^{i+1}\times\M^{e-i}/\M^{e+1-i}\to\M^e\subseteq\Soc(A)
$$
induced by the multiplication are non--degenerate.

In particular if $A$ is Gorenstein then $\dim_k(\M^i/\M^{i+1})=\dim_k(\M^{e-i}/\M^{e+1-i})$ for $i=0,\dots,e$.
\qed
\endproclaim

It follows that up to a proper choice of the variables we have an isomorphism $R\cong k[x_1,\dots,x_{d-2}]/I_{G_d}$ where
$$
I_{G_d}=(x_ix_j-\delta_{i,j}x_1^2)_{1\le i\le j\le d-2}\subseteq k[x_1,\dots,x_{d-2}],
$$
$\delta_{i,j}$ being the Kronecker symbol.

Since $I_{G_d}\subseteq k[x_1,\dots,x_{d-2}]\subseteq k[x_0,\dots,x_{d-2}]$ is a homogeneous ideal, then it defines a subscheme in $\p{d-2}$ with Hilbert function $h(0)=1$, $h(1)=d-1$ and $h(t)=d$ for $t\ge2$. It is obviously arithmetically Gorenstein due to the discussion above, hence $G_d$ defines a point in $\Hilb_{d}^{aG,0}(\p{d-2})$. 

\definition{Definition 3.2}
A $G$--fat point of degree $d$ in the projective space $\p{d-2}$ is any scheme projectively equivalent to the above scheme $G_d\subseteq\p{d-2}$.

We denote by ${\Cal G}_d\subset\Hilb_{d}^{aG}(\p{d-2})$ the subset of $G$--fat points. 
\enddefinition

The $G$ of $G$--fat stands for Gorenstein in order to distinguish the above scheme by the standard fat point in $\p{d-2}$  defined by the ideal $(x_1,\dots,x_{d-2})^2\subseteq k[x_0,\dots,x_{d-2}]$ (which has, by the way, degree $d-1$, not $d$). 

Now we come back to the homogeneous ideal $I_X\subseteq k[x_0,\dots,x_{d-2}]$ of our scheme $X$. Notice that by construction $I\subseteq I_{G_d}+(x_0)$. Since $x_0$ is regular and $I_X$ is minimally generated by $d(d-3)/2$ quadratic polynomials by sequence (2.3.1) above we finally obtain the first part of the statement of the following

\proclaim{Proposition 3.3}
Let $X\in \Hilb_{d}^{aG,0}(\p{d-2})$. Then there exist homogeneous coordinates $x_0,\dots,x_{d-2}$ in $\p{d-2}$ such that the homogeneous ideal $I_X\subseteq k[x_0,\dots,x_{d-2}]$ of $X$ is
$$
I_X=(x_ix_j-\delta_{i,j}x_1^2+x_0\ell_{i,j})_{1\le i\le j\le {d-2},\ (i,j)\ne(1,1)}
$$
where $\ell_{i,j}\in k[x_0,\dots,x_{d-2}]$ are linear forms for $i,j=1,\dots,{d-2}$.

Moreover if $X$ is irreducible one can also assume that $\ell_{i,j}\in k[x_1,\dots,x_{d-2}]$.
\endproclaim
\demo{Proof}
It remains to prove the statement in the irreducible case. In this case $X$ contains a unique point that we can assume to be $[1,0\dots,0]$ thus $I_X\subseteq (x_1,\dots,x_{d-2})$ and $x_0$ is a regular element in $S_X$.
\qed
\enddemo

\proclaim{Corollary 3.4}
${\Cal G}_d$ is contained in the intersection of all the irreducible components of $\Hilb_{d}^{aG,0}(\p{d-2})$. In particular each $X\in \Hilb_{d}^{aG,0}(\p{d-2})$ can be deformed to a point in ${\Cal G}_d$
\qed
\endproclaim

Trivially, if $\Hilb_{d}^{aG,0}(\p{d-2})$ is reducible, then every point in ${\Cal G}_d$ turns out to be obstructed. This is not evident if $\Hilb_{d}^{aG,0}(\p{d-2})$ is irreducible. 

 \proclaim{Theorem 3.5}
The dimension of the tangent space of $\Hilb_{d}^{aG,0}(\p{d-2})$ at the points in  ${\Cal G}_d$ is $(d^3-7d)/6$. Thus the points in  ${\Cal G}_d$ are obstructed if and only if $d\ge6$.
\endproclaim
\demo{Proof}
Assume that $\Hilb_{d}^{aG,0}(\p{d-2})$ is irreducible (this is known if $d=3,4,5$ and we will prove this later on in the case $d=6$). Hence $\Hilb_{d}^{aG,0}(\p{d-2})=\Hilb_{d}^{gen}(\p{d-2})$ which has dimension $d(d-2)$ (see Proposition 2.5 above). Since all the points $X\in{\Cal G}_d$ are projectively equivalent to $G_d$ it suffices to compute the dimension of the tangent space to $\Hilb_{d}(\p{d-2})$ at $G_d$, i.e. $h^0\big(\p{d-2},(\Im/\Im^2)\check{\ }\big)$, $\Im$ being the sheaf of ideals of $G_d$ as subscheme of $\p{d-2}$.

To this purpose consider the affine space $\a{d-2}=\p{d-2}\setminus\{\ x_0=0\ \}$. Since obstructedness does not depend on the chosen projective model of $X$ we can assume $X=G_d\subseteq\a{d-2}\subseteq\p{d-2}$, then $H^0\big(\p{d-2},(\Im/\Im^2)\check{\ }\big)=H^0\big(\a{d-2},(\Im/\Im^2)\check{\ }_{\vert\a{d-2}}\big)$ (see [Ha2], exercise III.2.3). Let $A:=k[x_1,\dots,x_{d-2}]$, so that  $G_d=\spec(A/I_{G_d})\subseteq\a{d-2}$, hence
$$
H^0\big(\p{d-2},(\Im/\Im^2)\check{\ }\big)\cong H^0\big(\a{d-2},\widetilde{(I_{G_d}/I_{G_d}^2)\check{\ }}\big)\cong\hom_{A/I_{G_d}}\big(I_{G_d}/I_{G_d}^2,A/I_{G_d}\big).
$$
Since $A/I_{G_d}$ is an Artinian Gorenstein algebra, then $\hom_{A/I_{G_d}}\big(I_{G_d}/I_{G_d}^2,A/I_{G_d}\big)$ and  $I_{G_d}/I_{G_d}^2$ have the same length (see [No], Chapter 5, Theorem 21), thus $h^0\big(\p{d-2},(\Im/\Im^2)\check{\ }\big)=\dim_k(I_{G_d}/I_{G_d}^2)$.

Since $I_{G_d}$ is generated by homogeneous polynomials of degree $2$ it is clear that $I_{G_d}/I_{G_d}^2$ contains the components of degree $2$ and $3$ of $I_{G_d}$. By its very definition in degree $3$ all the monomials in the variables $x_1,\dots,x_{d-2}$ are in $I_{G_d}$. On the other hand $I_{G_d}^2$ is generated by all the products of the generators of $I_{G_d}$. In particular it contains all the monomial containing at most a square, which come from the products of two products of distinct variables. Since $x_i^3x_j-x_1^2x_ix_j\in I_{G_d}^2$, $i=1,\dots,n$, it follows that $I_{G_d}^2$ also contains all the other monomials of degree $4$. We conclude that 
$$
\dim(I_{G_d}/I_{G_d}^2)={d\choose 3}+{{d-1}\choose2}-1={{d^3-7d}\over6}.
$$
Thus $G_d$ is obstructed if and only if $h^0\big(\p{d-2},{\Cal N}\big)=\dim_k(I_{G_d}/I_{G_d}^2)>d(d-2)$, i.e. if and only if $d\ge6$.
\qed
\enddemo

\head
4. Classification of degree $d\le6$ local Artinian Gorenstein $k$--algebras
\endhead
Consider the natural action of $\PGL_{d-1}$ on $\Hilb_{d}^{aG,0}(\p{d-2})$. Let $X',X''\in \Hilb_{d}^{aG}(\p{d-2})$: due to Corollary 2.3, if $X'\cong X''$ abstractly then $X',X''$ are projectively isomorphic as subschemes of $\p{d-2}$, thus the action is transitive and there exists a natural stratification
$$
\Hilb_{d}^{aG,0}(\p{d-2})=\bigcup_{A} O(A)
$$
where $A$ runs in the set of all possible non--isomorphic Artinian Gorenstein $k$--algebras of degree $d$ and $O(A)$ denotes the $\PGL_{d-1}$--orbit of any arithmetically Gorenstein embedding $\spec(A)\subseteq\p{d-2}$: e.g.
$$
{\Cal G}_d=O(k[x_0,\dots,x_{d-2}]/(x_ix_j-\delta_{i,j}x_1^2)_{\vert\ 1\le i\le j\le d-2}).
$$
As explained above there always is a surjective morphism $\PGL_{d-1}\to O(A)$, thus $O(A)$ is always irreducible and open in its closure $\overline{O(A)}$ inside $\Hilb_{d}^{aG,0}(\p{d-2})$.

It is then interesting to understand the intrinsic structure of such algebras $A$'s.
Since $A$ is Artinian then it is direct sum of local Artinian Gorenstein $k$--algebras of degree at most $d$, thus it is natural to begin our inspection from these elementary bricks. 

So, from now on in this section, we will assume that $A$ is a local Artinian Gorenstein $k$--algebra of degree $d$ with maximal
ideal $\M$.

In general we have a filtration
$$
A\supset\M\supset\M^2\supset\dots\supset\M^e\supset\M^{e+1}=0 
$$
for some integer $e\ge1$, so that its associated graded algebra
$$
\gr(A):=\bigoplus_{i=0}^\infty\M^i/\M^{i+1}
$$
is a vector space over $k\cong A/\M$ of finite dimension $d=\dim_k(A)=\dim_k(\gr(A))=\sum_{i=0}^e\dim_k(\M^i/\M^{i+1})$.

\definition{Definition 4.1}
Let $k$, $A$, $\M$ be as above. If $\M^e\ne0$ and $\M^{e+1}=0$ we define the
level of $A$ as $\lev(A):=e$.

If $\lev(A)=e$ and $n_i:=\dim_k(\M^i/\M^{i+1})$, $0\le i\le e$, we define the Hilbert function of $A$ as the vector $H(A):=(n_0,\dots,n_e)\in{\Bbb N}^{e+1}$.
\enddefinition

In any case $n_0=1$. Recall that $A$ is said to be Gorenstein if and only if its
socle $\Soc(A):=0\colon{\frak M}$ has dimension
$1$ over the residue field $A/{\frak M}$, which in our case is $k$. 
If $e=\lev(A)$, since $0\ne\M^e\subseteq\Soc(A)$, then equality must
hold i.e.
$n_e=1$. 

If $A$ is also graded then one can use Proposition 3.1 above. If $A$ is not graded this assertion is no longer true, so the classification of
non--graded local and Artinian Gorenstein algebras is a little bit more difficult (at least
for the computations). 

If either $e=0$ (resp. $e=1$) then trivially $A\cong A_{0,1}:=k$ (resp. $A\cong A_{1,2}:=k[x_1]/(x_1^2)$), so we can assume $e\ge3$ in what follows. Consider
the case $H(A)=(1,n,1,\dots,1)$, so that $\lev(A)=d-n\ge2$. Let $a_1,\dots,a_n$ be a minimal set of generators of 
$\M$ (so that their classes are linearly independent in $\M/\M^2$). If $n\ge2$, we can always make the assumption
$a_1^2\not\in\M^3$. Indeed if $a_i^2\in \M^3$ for each $i=1,\dots,n$
then at least one of the mixed products $a_ia_j$, where $i,j=1,\dots,n$ with $i<j$, cannot lie in $\M^3$. E.g. assume that $a_1a_2\not\in\M^3$: thus  $(a_1+a_2)^2\not\in\M^3$, whence we can change $a_1$ with $a_1+a_2$ obtaining in particular
$\M^2=(a_1^2)$.

We have relations of the form $a_1a_i=\alpha_ia_1^2$, $\alpha_i\in A$. By changing
$a_i$ with $a_i-\alpha_ia_1$ we can take the generators in such a way that
$a_1a_i=0$. It follows $\M^i=(a_1^i)$ for each $i\ge2$.

If $a_ia_j\in\M^t\setminus\M^{t+1}$, where $i,j=2,\dots,n$ and $t\le\lev(A)=d-n$, then $a_ia_j=\beta a_1^t$ and $\beta$ is
invertible thus $a_1^{t+1}=\beta^{-1}a_1a_ia_j=0$, whence $\M^{t+1}=0$,
which implies $t=\lev(A)=d-n$. It follows that $a_ia_j=\beta_{i,j}a_1^{d-n}$ for
some $\beta_{i,j}\in k\subseteq A$, where $i,j=2,\dots,n$ and $\beta_{i,j}=\beta_{j,i}$.

Let $B:=(\beta_{i,j})_{2\le i,j\le n}$. Assume $\det(B)=0$: we can find a non--zero
$\gamma:=(\gamma_2,\dots,\gamma_n)\in k^{n-1}$ such that $B{}^t\gamma=0$. If e.g. $\gamma_2\ne0$ we can substitute $a_2$ with $\gamma_2a_2+\dots+\gamma_na_n$, hence we can assume $\beta_{i,2}=0$ for each $i=2,\dots,n$, thus $V:=\langle
a_2,a_1^{d-n}\rangle\subseteq \Soc(A)$. On the other hand $a_2,a_1^{d-n}$ are both non--zero and they must be linearly independent otherwise $a_2\in\M^{d-n}\subseteq\M^2$ which is absurd. We conclude that $\dim_k(V)=2$ which is again an absurd since $\dim_k(\Soc(A))=1$. It follows that
$\det(B)\ne0$.

Set ${\frak N}:=(a_2,\dots,a_n)\subseteq A$. Notice that the classes of
$a_2,\dots a_n$ in $U:={\frak N}/({\frak N}\cap\M^2)$ form a basis of $U$. The
condition $\det(B)\ne0$ implies that the symmetric bilinear form 
$\psi\colon U \times U\to\M^{d-n}=\Soc(A)\cong k$ induced by the multiplication is
non degenerate. It then follows that we can change the generators $a_2,\dots,a_n$ of
$\frak N$ with suitable linear combinations in such a way that $a_ia_j=0$ (when
$1\le i\ne j\le n$) and $a_h^2=a_1^{d-n}$ (when $h=2,\dots,n$).

Now let $\varphi\colon k[x_1,\dots,x_n]\twoheadrightarrow A$ be the epimorphism defined by $x_i\mapsto a_i$. Then $I:=(x_ix_j-\delta_{i,j}x_1^{d-n},x_1^{d-n+1})_{i,j=1,\dots,n\ (i,j)\ne(1,1)}\subseteq\ker(\varphi)$ and 
$\dim_k(k[x_1,\dots,x_n]/I)=d$ thus the induced map $\overline{\varphi}\colon
k[x_1,\dots,x_n]/I\twoheadrightarrow A$ must be an isomorphism.

\proclaim{Theorem 4.2}
If $H(A)=(1,n,1,\dots,1)$, $1\le n\le d-1$, then
$$
A\cong A_{n,d}:=k[x_1,\dots,x_n]/(x_ix_j-\delta_{i,j}x_1^{d-n},x_1^{d-n+1})_{1\le i\le j\le {n},\ (i,j)\ne(1,1)}.
$$
\endproclaim
\demo{Proof}
It remains only to examine the case $n=1$: this is equivalent to $\lev(A)=d-1$ and it is not difficult to prove that
$A\cong k[x]/(x^d)$ which is $A_{1,d}$.
\qed
\enddemo

\remark{Remark 4.3}
The two extremal cases are interesting, namely $A_{1,d}\cong k[x]/(x^d)$ and $A_{d-2,d}\cong
k[x_1,\dots,x_{d-2}]/(x_ix_j,x_i^2-x_1^2)_{1\le i\le j\le {d-2}}$.

In both these cases $A$ turns out to be naturally graded. In the remaining
cases $A$ is obviously non--graded.
\endremark
\medbreak

If $d=3,4,5$ then $A$ is completely described
by Theorem 4.2 above. It is then natural to ask the following

\definition{Question 4.4}
Is it possible to classify Artinian Gorenstein $k$--algebras up to isomorphism when $d\ge6$?
\enddefinition

In the remaining part of this section we will answer the above question in the first case $d=6$. We do not know any general answer for greater values of $d$. We must then have $e:=lev(A)\le6$, hence it is easy to check that if $\lev(A)\ne 3$ then $H(A)=(1,1,1,1,1,1)$ (if $\lev(A)=5$), $H(A)=(1,2,1,1,1)$ (if $\lev(A)=4$), $H(A)=(1,4,1)$ (if $\lev(A)=2$): all these cases can be handled with the help of Theorem 4.2. When $\lev(A)=3$ then either $H(A)=(1,3,1,1)$ and again we can treat this case using Theorem 4.2, or $H(A)=(1,2,2,1)$.

In this case we have
$\M=(a,b)$, hence $\M^2=(a^2,ab,b^2)$. Since $\dim_k(\M^2/\M^{3})=2$ there is a non--trivial relation
$\alpha a^2+\beta  b^2+\gamma ab\in\M^3$, where as usual $\alpha,\beta,\gamma\in
k\subseteq A$. Notice that we can assume $\gamma\ne0$: if not and e.g. $\beta\ne0$ we change $b$ with $a+b$ obtaining a relation as above with $\gamma\ne0$.
Since $\gamma\ne0$ we can make the position $\gamma=2$ from now on. We distinguish two cases according that $\alpha \beta =1$ or $\alpha\beta\ne1$.

In the first case, changing $b$ with $\alpha a+\beta b$, we can assume $b^2\in\M^3$, thus $\M^2=(a^2,ab)$,
$\M^3=(a^3,a^2b)$: Nakayama's lemma yields the existence of $\delta,\varepsilon\in k$, not both zero, such that
$\delta a^3+\varepsilon a^2b\in\M^4=0$.

If $\varepsilon\ne0$, then we also have $\M^3=(a^3)$. If $\varepsilon=0$ then $a^3=0$ and $\M^3=(a^2b)$. Notice that $(a+b)^3=a^3+3a^2b+3ab^2+b^3$: since $a^3=0$ and $ab^2,b^3\in\M^4=0$ then $(a+b)^3=3a^2b\in\M^3\setminus\M^4$. Changing $a$ with $a+b$ it is then easy to check that $a^2$ and $ab$ are again linearly independent in $\M^2$ modulo $\M^3$, Thus again $\M=(a,b)$, $\M^2=(a^2,ab)$, $\M^3=(a^3)$.

Since in both the cases we have $b^2,a^2b\in (a^3)$ we obtain an isomorphism $k[x_1,x_2]/I\cong A$ where $I:=(x_2^2-\eta x_1^3, x_1^2x_2-\vartheta x_1^3, x_1^4, x_2^4)$, where $\eta,\vartheta\in k$ since ${\frak M}^4=0$. One has $\dim_k(\Soc(k[x_1,x_2]/I))=1$ if and only if $\vartheta\ne0$, so we can take
$\vartheta=1$. The sequence of transformations $(x_1,x_2)\mapsto \eta^{-1}(x_1,x_2)$ followed by $(x_1,x_2)\mapsto(x_1,x_2+x_1^2/2)$ allows us to assume that $\eta=0$. Now the transformation $(x_1,x_2)\mapsto (x_1-x_2/3,x_2)$ finally yields
$$
A\cong A_1^{sp}:=k[x_1,x_2]/(x_2^2,x_1^3).
$$

In the second case, i.e. $\alpha \beta \ne1$, let $\lambda\in k$ be a non--zero
root of
$P(t):=t^2-2t+\alpha \beta $ and set $a_1:=\alpha a+\lambda b$ and
$b_1=\lambda a+\beta  b$. On one hand there is $m\in\M^3$ such that $\alpha \lambda a^2+\beta \lambda b^2=-2\lambda ab+m$, whence
$$
a_1b_1=\alpha \lambda a^2+(\alpha \beta +\lambda^2)ab+\beta \lambda
b^2=abP(\lambda)+m=m\in\M^3.
$$
On the other hand
$$
D:=\left\vert\matrix
\alpha &\lambda\\
\lambda&\beta 
\endmatrix\right\vert=\alpha \beta -\lambda^2=2(\alpha \beta-\lambda) .
$$
$D\in\frak M$ if and only if $D=0$ (since $D\in k$), i.e. $\lambda^2=\alpha \beta =\lambda$, which
is not possible by the assumptions.

An immediate consequence of the above discussion is that we can substitute $a,b$ with $a_1,b_1$, hence $\M=(a,b)$, $\M^2=(a^2,b^2)$ (since $ab\in\M^3$ by the choice of $\lambda$) and $\M^3=(a^3)$
(exchanging possibly the roles of $a$ and $b$). Moreover we have the relations
$ab=\mu a^3$, $b^3=\nu a^3$, $a^4=0$ where as above we take
$\mu,\nu\in k\subseteq A$.

We conclude that there is 
$\varphi\colon k[x_1,x_2]\twoheadrightarrow A$ given by $x_1\mapsto a$, $x_2\mapsto b$
and satisfying
$I=(x_1x_2-\mu x_1^3,x_2^3-\nu x_1^3,x_1^4)\subseteq\ker\varphi$.
The induced map $\overline{\varphi}\colon k[x_1,x_2]/I\to A$ must be an isomorphism
by degree computation. Moreover $\dim_k(Soc(k[x_1,x_2]/I)=1$ if and only if $\nu\ne0$: changing suitably $x_1$ we can thus assume $\nu=1$. The sequence of transformations $(x_1,x_2)\mapsto \mu^{-1}(x_1,x_2)$ followed by $(x_1,x_2)\mapsto (x_1,x_2+x_1^2)$ allows us to assume that $\mu=0$. Finally the transformation $(x_1,x_2)\mapsto(x_1-x_2,x_1+x_2)$ yields
$$
A\cong A_2^{sp}:=k[x_1,x_2]/(x_2^2-x_1^2,x_1^3).
$$

\proclaim{Proposition 4.5}
$A_1^{sp}\not\cong A_2^{sp}$.
\endproclaim
\demo{Proof}
Notice that in $A_1^{sp}$ there is $\overline{x}_2\in\M\setminus \M^2$ such that $\overline{x}_2^2=0$.
On the other hand in $A_2^{sp}$ let $y\in \M\setminus\M^2$ be such that $y^2=0$. If $y=\alpha \overline x_1+\beta\overline x_2$ necessarily $\alpha,\beta$ are not both in $\M$.

If e.g. $\alpha\not\in\M$, then $\alpha^2\not\in\M$ too. Then the condition
$$
(\alpha^2+\beta^2)\overline{x}_2^2+2\alpha\beta\overline{x}_1\overline{x}_2=
(\alpha\overline{x}_1+\beta\overline{x}_2)^2=0
$$
implies $\alpha^2+\beta^2,\alpha\beta\in\M$. We thus obtain $\beta\in\M$ (since $\alpha\not\in\M$) whence $\alpha^2\in\M$ an absurd.
\qed
\enddemo

\remark{Remark 4.6}
An immediate consequence of the description above is the existence of exactly $20$ non--isomorphic Artinian Gorenstein $k$--algebras of degree $6$, namely $A_{0,1}^{\oplus 6}$, $A_{0,1}^{\oplus 4}\oplus A_{1,2}$, $A_{0,1}^{\oplus 3}\oplus A_{1,3}$, $A_{0,1}^{\oplus 2}\oplus A_{1,2}^{\oplus2}$, $A_{0,1}^{\oplus 2}\oplus A_{1,4}$, $A_{0,1}^{\oplus 2}\oplus A_{2,4}$, $A_{0,1}\oplus A_{1,2}\oplus A_{1,3}$, $A_{1,2}^{\oplus3}$, $A_{0,1}\oplus A_{1,5}$, $A_{0,1}\oplus A_{2,5}$, $A_{0,1}\oplus A_{3,5}$, $A_{1,2}\oplus A_{1,4}$, $A_{1,2}\oplus A_{2,4}$, $A_{1,3}^{\oplus 2}$, $A_{1,6}$, $A_{2,6}$, $A_{3,6}$, $A_{4,6}$, $A_1^{sp}$, $A_2^{sp}$.
\endremark
\medbreak

\remark{Remark 4.7}
When $d\ge7$ we have many new possible shapes for the Hilbert function $H(A)$, but not all of them can really occur. For example if $d=7$ other the cases $A_{n,7}$ and $H(A)=(1,2,2,1,1)$, which can be described as the $A_i^{sp}$, we can also have $H(A)=(1,3,2,1)$ and $H(A)=(1,2,3,1)$, but the second one cannot occur in the Gorenstein case. Indeed all the relations are in degree $3$, thus they are all homogeneous. In particular such an algebra should be graded, hence $H(A)$ should be symmetric by Proposition 3.1.
\endremark
\medbreak

\head
5. Irreducibility of $\Hilb_{6}^{aG}(\p{4})$ and its singular locus
\endhead

Since there are no general results about $\Hilb_{d}^{aG}(\p{d-2})\subseteq \Hilb_{d}(\p{d-2})$ when $d\ge6$, it is quite natural to study the first case, namely $\Hilb_{6}^{aG}(\p{4})\subseteq \Hilb_{6}(\p{4})$.

In Section 2 we defined $\Hilb_{6}^{gen}(\p{4})\subseteq\Hilb_{6}^{aG}(\p{4})$ as follows.
Let ${\Cal W}\subseteq\Sim^6(\p{4})$ be the irreducible open subset parametrizing cycles of $6$ distinct points. The natural morphism $s\colon \Hilb_{6}(\p{4})\to\Sim^6(\p{4})$ induces an isomorphism $s^{-1}({\Cal W})\to{\Cal W}$. Then we defined the non--empty subscheme $\Hilb_{6}^{gen}(\p{4}):=\Hilb_{6}^{aG}(\p{4})\cap \overline{s^{-1}({\Cal W})}$. 

\proclaim{Proposition 5.1}
$\Hilb_{6}^{aG}(\p{4})=\Hilb_{6}^{gen}(\p{4})$ hence it is irreducible. $\Hilb_{6}^{aG,0}(\p{4})$ is open and dense inside $\Hilb_{6}^{gen}(\p{4})$.
\endproclaim
\demo{Proof}
We will then prove that each $X\in \Hilb_{6}^{aG}(\p{4})$ is specialization of a family ${\Cal X}\to B\subseteq\a1$ of points of $\Hilb_{6}^{gen}(\p{4})$, thus $\Hilb_{6}^{aG}(\p{4})=\Hilb_{6}^{gen}(\p{4})$ is irreducible since $\Hilb_{6}^{gen}(\p{4})$ is closed in $\Hilb_{6}^{aG}(\p{4})$ by definition. To this purpose we first notice that $\Hilb_{6}^{aG}(\p{4})=\Hilb_{6}^{aG,0}(\p{4})\cup\Hilb_{6}^{aG,2}(\p{4})\cup\Hilb_{6}^{aG,3}(\p{4})$ due  Proposition 2.1. 

If $X\in \Hilb_{6}^{aG,2}(\p{4})$ then it is the complete intersection of a conic and a cubic in a plane. If $X\in \Hilb_{6}^{aG,3}(\p{4})$ then it is a divisor of degree $6$ on a line. It is clear that smooth schemes of these kinds are in $\Hilb_{6}^{gen}(\p{4})$: moreover it is not difficult to check that any other scheme of these kind can be obtained as specialization of the smooth ones. An easy parameters computation shows that $\dim(\Hilb_{6}^{aG,2}(\p{4}))=17$, $\dim(\Hilb_{6}^{aG,3}(\p{4}))=12$.

We now examine points $X\in \Hilb_{6}^{aG,0}(\p{4})$, which is open inside $\Hilb_{6}^{gen}(\p{4})$ by Proposition 2.5 and the definition of $\Hilb_{6}^{gen}(\p{4})$. To this purpose we will make use of Corollary 2.4, in the following way. Let $X\in\Hilb_{6}^{aG,0}(\p{4})$. Then $X\cong \spec(A)$ for a suitable Artinian Gorenstein $k$--algebra of degree $6$. Making use of the classification given in the previous section, we will construct a family ${\Cal A}\to\a1$ of Artinian Gorenstein degree $6$ quotients of a fixed polynomial ring $k[b,x_1,\dots,x_n]$, whose fibre over $b=0$ is $A$ and over $b\in\a1\setminus\{\ 0\ \}$ corresponds to a point in $\Hilb_{6}^{gen}(\p{4})$. Then, if $B:=\a1$, $\spec({\Cal A})\subseteq\a n\times B$ satisfies Corollary 2.4. In particular we have a flat family ${\Cal X}'$ over an open neighborhood of $b=0$ in $\a1$ whose special fibre is projectively isomorphic in $\p4$ to $X$ and the other fibres represent points in $\Hilb_{6}^{gen}(\p{4})$. Since $\Hilb_{6}^{gen}(\p{4})$ is the closure of $O(A_{0,1}^{\oplus 6})$ then it contains the whole of $O(A)$, hence $X$.

Assume that the tangent spaces at the points of $X=\spec(A)$ have all dimension at most $1$. Then $A$ is one of the following: $A_{0,1}^{\oplus 6}$, $A_{1,2}^{\oplus e}\oplus A_{0,1}^{\oplus 6-2e}$ for $e=1,2,3$,  $A_{1,3}^{\oplus e}\oplus A_{0,1}^{\oplus 6-3e}$ for $e=1,2$,  $A_{1,3}\oplus A_{1,2}\oplus A_{0,1}$, $A_{1,4}\oplus A_{0,1}^{\oplus2}$, $A_{1,4}\oplus A_{1,2}$, $A_{1,5}\oplus A_{0,1}$, $A_{1,6}$. By its very definition the general point of $\Hilb_{6}^{gen}(\p{4})$ corresponds to $A_{0,1}^{\oplus 6}$: notice that $A_{0,1}^{\oplus 6}\cong k[x]/(p(x))$ for some polynomial $p(x):=\prod_{h=0}^5(x-\alpha_h)$ without multiple roots.

In the other cases it suffices to show the existence of an ${\Cal A}\to\a1$ with general fibre $A_{0,1}^{\oplus 6}$ and special fibre the fixed $A$. We examine only the case $A\cong A_{1,6}\cong k[x]/(x^6)$ the other ones being similar. In this case consider the family over $\a1$ with coordinate $b$ given by ${\Cal A}:=k[b,x]/(\prod_{h=0}^5(x-b\alpha_h))$.

Now assume that there exists a point in $X$ with tangent space of dimension at least $2$. Then $A$ is one of the following: $A_{2,4}\oplus A_{0,1}^{\oplus2}$, $A_{2,4}\oplus A_{1,2}$, $A_{2,5}\oplus A_{0,1}$, $A_{2,6}$, $A_1^{sp}$, $A_2^{sp}$, $A_{3,5}\oplus A_{0,1}$, $A_{3,6}$, $A_{4,6}$. We examine again only the first case. Let $I:=(x_2^2-x_1^2-x_2^4+bx_1^6,x_1^3-bx_1^7,x_1x_2,bx_2)$ and consider the family ${\Cal A}:=k[b,x_1,x_2]/I$. We have 
$$
\gather
{\Cal A}_0=k[x_1,x_2]/((x_2^2-x_1^2,x_1x_2)\cap(x_1,x_2-1)\cap (x_1,x_2+1))\cong A_{2,4}\oplus A_{0,1}^{\oplus2},\\
{\Cal A}_b=k[x_1,x_2]/(-x_1^2+bx_1^6,x_2)\cong A_{1,2}\oplus A_{0,1}^{\oplus4},\quad \text{if $b\ne0$}.
\endgather
$$
Since we have just proved that $A_{1,2}\oplus A_{0,1}^{\oplus4}$ defines a point in $\Hilb_{6}^{gen}(\p{4})$ we have completed the proof also in this case.
\qed
\enddemo

Let $A$ be an Artinian Gorenstein $k$--algebra of degree $6$. We already proved in Theorem 3.5 that if $X\in{\Cal G}_6$ then it is obstructed. The aim of the remaining part of the present section is to show that in all the other cases $X$ is non--obstructed. To this purpose we have to compute the dimension of the tangent space to $\Hilb_{6}^{aG}(\p{4})$ at $X$, i.e. $h^0\big(\p{d-2},(\Im/\Im^2)\check{\ }\big)$, $\Im$ being the sheaf of ideals of $X$ as subscheme of $\p{d-2}$ (see the proof of Theorem 3.5).

We begin by examining the case of irreducible schemes in $\Hilb_{6}^{aG,0}(\p{4})$. Our first step is to compute their projective models in $\p4$.

\proclaim{Lemma 5.2}
Let $X\in\Hilb_{6}^{aG,0}(\p{4})$ be  irreducible. Then, up to a proper choice of the coordinates $x_0,x_1,x_2,x_3,x_4$ in $\p{4}$, the homogeneous ideal of $X\subseteq\p4$ is:
\item{i)} $(x_1x_2-x_0x_3,x_1x_3-x_0x_4,x_1x_4-x_2x_3,x_2x_4,x_3x_4,x_1^2-x_0x_2,x_2^2-x_0x_4,x_3^2,x_4^2)$ if $X\in O(A_{1,6})$;
\item{ii)}$(x_1x_2-x_0x_3,x_1x_3-x_4^2,x_1x_4,x_2x_3,x_2x_4,x_3x_4,x_1^2-x_0x_2,x_2^2-x_4^2,x_3^2)$ if $X\in O(A_{2,6})$;
\item{iii)} $(x_1x_2-x_4^2,x_1x_3,x_1x_4,x_2x_3,x_2x_4,x_3x_4,x_1^2-x_0x_2,x_2^2,x_3^2-x_4^2)$ if $X\in O(A_{3,6})$;
\item{iv)} $(x_1x_2,x_1x_3,x_1x_4,x_2x_3,x_2x_4,x_3x_4,x_2^2-x_1^2,x_3^2-x_1^2,x_4^2-x_1^2)$ if $X\in O(A_{4,6})$;
\item{v)} $(x_1x_2,x_1x_3-x_0x_4,x_1x_4-x_2x_3,x_2x_4,x_3x_4,x_1^2-x_0x_2,x_2^2,x_3^2,x_4^2)$ if $X\in O(A_1^{sp})$;
\item{vi)} $(x_1x_2,x_1x_3-x_0x_4,x_1x_4-x_2x_3,x_2x_4,x_3x_4,x_1^2-x_0x_2,x_2^2,x_3^2-x_0x_2,x_4^2)$ if $X\in O(A_2^{sp})$.
\endproclaim
\demo{Proof}
We have to show that the homogeneous coordinate ring $S_X$ is Gorenstein in each of the above cases. Due to exercise 3.6.22 of [B--H] it suffices to show that $(S_X)_{(x_0,\dots,x_4)}$ is Gorenstein. Due to Proposition 3.1.19 (b) of [B--H] it suffices to find a regular element $x\in (S_X)_{(x_0,\dots,x_4)}$ such that  $(S_X/(x))_{(x_0,\dots,x_4)}$ is Gorenstein. 

We describe only the case i), the other ones being similar. First of all notice that $X$ is concentrated at the point $[1,0,0,0,0]$. By dehomogenizing with respect to $x_0$ we obtain 
$$
X\cong\spec(k[x_1,\dots,x_4]/(x_1^2-x_2,x_1^3-x_3,x_1^4-x_4,x_1^6))\cong\spec(A_{1,6}).
$$

Since $x_0$ is a regular element in $S_X$ it suffices to check that $(S_X)_{(x_0,\dots,x_4)}$ is Gorenstein. Since
$$
S_X/(x_0)\cong k[x_1,\dots,x_4]/(x_1x_2,x_1x_3,x_1x_4-x_2x_3,x_2x_4,x_3x_4,x_1^2,x_2^2,x_3^2,x_4^2)
$$
is local, then it suffices to verify that $S_X/(x_0)\cong (S_X)_{(x_0,\dots,x_4)}$ is Gorenstein. It is easy to check that $\lev(S_X/(x_0))=2$ and that its maximal ideal is ${\frak M}=(x_1,x_2,x_3,x_4)$ and ${\frak M}^2=(x_2x_3)$. Now since the matrix of the bilinear form on ${\frak M}/{\frak M}^2\times{\frak M}/{\frak M}^2$ induced by the multiplication in $S_X/(x_0)$ is
$$
\pmatrix
0&0&0&1\\
0&0&1&0\\
0&1&0&0\\
1&0&0&0
\endpmatrix
$$
Proposition 3.1 implies that $S_X/(x_0)$, hence $S_X$, is Gorenstein.
\qed
\enddemo

We begin our study of the singular locus $\Sing(\Hilb_{6}^{aG}(\p{4}))$ of $\Hilb_{6}^{aG}(\p{4})$ by proving the following

\proclaim{Lemma 5.3}
An irreducible $X\in\Hilb_{6}^{aG,0}(\p{4})$ is obstructed if and only if $X\in{\Cal G}_6$.
\endproclaim
\demo{Proof}
Let $A$ be a local Artinian Gorenstein algebra, $X\in O(A)\subseteq\Hilb_{6}^{aG,0}(\p{4})$. Then we have to prove that $X$ is obstructed if and only if $A\cong A_{4,6}$. The \lq\lq if\rq\rq\ part is Theorem 3.5. In order to prove the \lq\lq only if\rq\rq\ part we distinguish the two cases $X\in O(A_{1,6})\cup O(A_{2,6})\cup O(A_{1}^{sp})\cup O(A_{2}^{sp})$ and $X\in O(A_{3,6})$, proving that
$$
h^0\big(\p{d-2},(\Im/\Im^2)\check{\ }\big)\le24=\dim(\Hilb_{6}^{aG}(\p{4})).
$$
Taking into account Lemma 5.2 this can be done by using the computer algebra software Macaulay [B--S]. In order to make the proof readable also for those ones who are not expert in Macaulay we produce here a more explicit argument.

As in the proof of Theorem 3.5 we can assume $X\cong \proj(k[x_0,\dots,x_4]/I)$ where $I\subseteq k[x_0,\dots,x_4]$ is the homogeneous ideal of $X$ described in Lemma 5.2, since obstructedness does not depend on the chosen projective model of $X$. In all the cases $X$ is then concentrated at the point $[1,0,0,0,0]$. Consider the affine space $\a{4}=\p{4}\setminus\{\ x_0=0\ \}$. Since $X\subseteq\a{4}\subseteq\p{4}$, then $H^0\big(\p{4},(\Im/\Im^2)\check{\ }\big)=H^0\big(\a{4},(\Im/\Im^2)\check{\ }_{\vert\a{4}}\big)$ (see [Ha2], Exercise III.2.3). We have $X=\spec(A/I)\subseteq\a{4}$, for a suitable ideal $I\subseteq A:=k[x_1,\dots,x_4]$ and we have only to compute $\dim_k(I/I^2)$ (see the analogous claim in the proof of Theorem 3.5).

Let $X\in O(A_{1,6})\cup O(A_{2,6})\cup O(A_{1}^{sp})\cup O(A_{2}^{sp})$. Then
$$
\gather
I=(x_1^2-x_2,x_1^3-x_3,x_1^4-x_4,x_1^6)\qquad \text{if $X\in O(A_{1,6})$,}\\
I=(x_1^2-x_2,x_1^3-x_3,x_1^4-x_2^2,x_1x_4)\qquad \text{if $X\in O(A_{2,6})$,}\\
I=(x_1^2-x_2,x_1x_3-x_4,x_3^2,x_1^3)\qquad \text{if $X\in O(A_{1}^{sp})$,}\\
I=(x_1^2-x_2,x_1x_3-x_4,x_1^2-x_3^2,x_1^3)\qquad \text{if $X\in O(A_{2}^{sp})$.}
\endgather
$$
It turns out that in the considered cases $X\subseteq\a{4}$ is a complete intersection subscheme, thus $I/I^2\cong A^{\oplus4}$. It then follows that $h^0\big(\p{d-2},(\Im_0/\Im_0^2)\check{\ }\big)=\dim_k(A^{\oplus4})=24$.

Let $X\in O(A_{3,6})$. By dehomogenizing with respect to $x_0$ we get 
$$
I=(x_1^2-x_2,x_1x_3,x_1x_4,x_3x_4,x_3^2-x_4^2,x_1x_2-x_4^2)
$$
and $\dim_k(I/I^2)=\dim_k(A/I^2)-\dim_k(A/I)=\dim_k(A/I^2)-6$: thus $X$ is unobstructed if $\dim_k(A/I^2)\le30$. First of all we will show that the set $C$ of the classes of the monomial of degree at most $3$ generate $A/I^2$. To this purpose we begin by noticing that since
$$
x_1^4-2x_1^2x_2+x_2^2,
x_1^3x_2-x_1^2x_4^2-x_1x_2^2+x_2x_4^2,x_1^2x_4^2,
x_1^3x_3-x_1x_2x_3,x_1^3x_4-x_1x_2x_4\in I^2
$$
then
$x_1^4,x_1^3x_2,x_1^3x_3,x_1^3x_4\in (2x_1^2x_2-x_2^2,x_1x_2^2-x_2x_4^2,x_1x_2x_3,x_1x_2x_4)+{I^2}$, all the remaining monomials of degree $4$ being in $I^2$ (this can be checked by hands for all the other $31$ monomial). Analogously $x_1^5\in (3x_1x_2^2-2x_2x_4^2)+{I^2}$, all the remaining monomials of degree $5$ being in $I^2$. Finally all the monomials of degree at least $6$ are in $I^2$, thus all the elements in $A/I^2$ can be represented by a polynomial of degree at most $3$. 
Since
$$
x_1^2x_3x_4-x_2x_3x_4, x_1^2x_3^2-x_1^2x_4^2-x_2x_3^2+x_2x_4^2\in I^2,
$$
then $x_2x_3x_4,-x_2x_3^2+x_2x_4^2\in{I^2}$.
Finally since
$$
x_2^2x_i-2x_1^2x_2x_i+x_1^4x_i\in I^2,\qquad i=2,3,4,
$$
then $x_2^2x_i\in (x_1x_2^2-x_2x_4^2,x_1x_2x_3,x_1x_2x_4)+{I^2}$, for $i=2,3,4$.
It follows that $A/I^2$ is generated by all the monomials in $A$ of degree at most $3$ but $x_2x_3x_4,x_2x_3^2,x_2^3,x_2^2x_3,x_2^2x_4$, whence $\dim_k(A/I^2)\le30$.
\qed
\enddemo

We now want to examine the general case. To this purpose we will make use of Corollaries 2.3 and 2.4 following the same strategy used in the proof of Proposition 5.1 (and Lemma 5.3), that is we will construct a flat family ${\Cal A}\to\a1$ whose general fibre is $A$ and special fibre is one of the algebras $A_{1,6}$, $A_{2,6}$, $A_{3,6}$, $A_{1}^{sp}$, $A_{2}^{sp}$. Since the special fibre of this family is non--obstructed in $\Hilb_{6}^{gen}(\p{4})$ the same is true for the general one. We thus obtain the following

\proclaim{Theorem 5.4}
$\Sing(\Hilb_{6}^{aG}(\p{4}))={\Cal G}_6$.
\endproclaim
\demo{Proof}
Again if $X\in \Hilb_{6}^{aG,2}(\p{4})\cup\Hilb_{6}^{aG,3}(\p{4})$, then $X\subseteq\a4$ is a complete intersection subscheme (see the proof of Proposition 5.1), thus we argue as in the proof of the case $X\in O(A_{1,6})\cup O(A_{2,6})\cup O(A_{1}^{sp})\cup O(A_{2}^{sp})$ of the above Lemma.

When $X\in \Hilb_{6}^{aG,0}(\p{4})$ we have already explained the idea of the proof. It remains to construct the family $\Cal A$ case by case. If $A\cong A_{0,1}^{\oplus 6}$ then trivially $A\cong k[x]/(\prod_{h=0}^5(x-\alpha_h))$. Then we consider the family over $\a1$ with coordinate $b$ given by ${\Cal A}:=k[b,x]/(\prod_{h=0}^5(x-b\alpha_h))$. In a similar way we can build families when $A$ is one of the following: $A_{1,2}^{\oplus e}\oplus A_{0,1}^{\oplus 6-2e}$ for $e=1,2,3$,  $A_{1,3}^{\oplus e}\oplus A_{0,1}^{\oplus 6-3e}$ for $e=1,2$,  $A_{1,3}\oplus A_{1,2}\oplus A_{0,1}$, $A_{1,4}\oplus A_{0,1}^{\oplus2}$, $A_{1,4}\oplus A_{1,2}$, $A_{1,5}\oplus A_{0,1}$.

If $A\cong A_{2,4}\oplus A_{0,1}^{\oplus2}$, since
$$
(x_1x_2,x_1^2b^2 - x_2^2,x_1^3)\cap(x_1^2+b^2,x_2)=( x_1x_2,x_1^4 + x_1^2b^2-x_2^2)
$$
then the family ${\Cal A}:= k[b,x_1,x_2]/( x_1x_2,x_1^4 + x_1^2b^2-x_2^2,x_2^3 )$ has general fibre $A$ and special fibre $A_{2,6}$. Similarly for $A_{2,4}\oplus A_{1,2}$. If $A\cong A_{2,5}\oplus A_{0,1}$, since
$$
(x_1x_2,x_1^3b - x_2^2,x_1^4)\cap(x_1-b,x_2)=(x_1x_2,x_2^2-bx_1^3+x_1^4)
$$
one can take ${\Cal A}:=k[b,x_1,x_2]/(x_1x_2,x_2^2-bx_1^3+x_1^4)$.  If $A\cong A_{3,5}\oplus A_{0,1}$, since
$$
\align
(x_1x_2,x_1x_3,x_2x_3,x_2^2-x_3^2,x_1^2b-x_3^2,x_1^3)&\cap(x_1+b,x_2,x_3)=\\
&=(x_1x_2,x_1x_3,x_2x_3,x_2^2-x_3^2,x_3^2-bx_1^2-x_1^3)
\endalign
$$
one can take ${\Cal A}:=k[b,x_1,x_2,x_3]/(x_1x_2,x_1x_3,x_2x_3,x_2^2-x_3^2,x_3^2-bx_1^2-x_1^3)$. 
\qed
\enddemo

\remark{Remark 5.5}
Each $X\in\Hilb_{d}^{aG,0}(\p{d-2})$ can be deformed to the the $G$--fat point by Corollary 3.4: as in the proofs above,  for each such $X$ one can easily obtain families ${\Cal X}\to\a 1$ whose general point is projectively equivalent to $X$ and whose special point is in ${\Cal G}_d$. Since each point in ${\Cal G}_d$ is smooth for $d=4,5$ by Theorem 3.5, it follows that the same is true for $X$. By examining also the points in $\Hilb_{d}^{aG,r}(\p{d-2})$, $r\ge1$, one obtains that $\Hilb_{d}^{aG}(\p{d-2})$ is smooth and irreducible. This result is a weaker form of the already mentioned irreducibility and smoothness of $\Hilb_4(\p{2})$ (see [A--V]) and it slightly extends the result proved in [MR] for Hilbert schemes of arithmetically Gorenstein of degree $d$ and codimension $3$ in $\p n$, $n\ge4$, to the case $n=3$ and $d=5$.
\endremark
\medbreak

\head
6. Example of constructions of schemes in $\Hilb_{6}^{aG,0}(\p{4})$
\endhead

In this section we present some examples of possible constructions of arithmetically Gorenstein schemes of dimension $0$ and degree $6$ in $\p4$, i.e. of points in $\Hilb_{6}^{aG,0}(\p{4})$. Moreover we will inspect which kind of schemes can be obtained with the different constructions we are dealing with.

\subhead
6.1. The Japanese construction
\endsubhead
Let $C$ and $F$ be respectively a conic and a cubic in $\p2$ without common components. Then $X=C\cap F$ is a scheme of degree $6$ and dimension $0$. Thus the Veronese embedding $v_2\colon \p2\to\p5$ embeds $X$ as a hyperplane section of the Veronese surface $V\subseteq\p5$, hence gives rise to a subscheme  in $\p4$, again denoted by $X$, of degree $6$ and dimension $0$: hence $X\subseteq\p4$ is arithmetically Cohen--Macaulay.  Since its canonical sheaf is $\omega_X\cong\Ofa{\p2}(2)$ then $\omega_X\cong\Ofa{\p4}(1)$, hence such a scheme is also subcanonical whence it is actually arithmetically Gorenstein.

\definition{Definition 6.1.1}
The construction described above is called Japanese.
\enddefinition

We call such construction \lq\lq Japanese\rq\rq\  since it is related to a suitable section of the Veronese surface, whose homogeneous ideal has a minimal free resolution described in [G--T].

The Japanese construction does not give all the points $X\in \Hilb_{6}^{aG,0}(\p{4})$ as proved in the following

\proclaim{Proposition 6.1.2}
$X\in \Hilb_{6}^{aG,0}(\p{4})$ can be obtained via the Japanese construction if and only if $X\not\in O(A_{2,6})\cup O(A_{3,5}\oplus A_{0,1})\cup O(A_{3,6})\cup O(A_{4,6})$.
\endproclaim
\demo{Proof}
All the tangent spaces at the points of such $X$ have dimensions at most $2$, thus $A\not\cong A_{3,5}\oplus A_{0,1}, A_{3,6}, A_{4,6}$. Moreover if $X$ is irreducible then $A$ must be graded: thus the vector $H(A)$, which does not depend on the grading, must be symmetric by Proposition 3.1, whence also $A_{2,6}$ cannot be obtained via the Japanese construction.

If all such tangent spaces have dimension at most one, then $A$ is direct sum of local algebras of type $A_{0,1}$ and $A_{1,d}$ then $X\subseteq\p1$. Thus the $2$--tuple embedding $\p1\to\p2$ embeds $X$ in $\p2$ as the subscheme of an integral conic $C$ and $X\in\left\vert \Ofa C(3)\right\vert$: since the restriction map $H^0\big(\p2,\Ofa{\p2}(3)\big)\to H^0\big(C,\Ofa C(3)\big)$ is surjective then there also exists a cubic $F$ such that $X=C\cap F$:

When at least one point in $X$ has tangent space of dimension $2$, then the homogeneous ideal of $X$ in $\p2$ coincides, up to projectivities, with $(x_1x_2,x_0x_1^2-x_0x_2^2)$, if $A\cong A_{2,4}\oplus A_{0,1}^{\oplus2}$, $(x_2^2,x_0x_1^2)$, if $A\cong A_{2,4}\oplus A_{1,2}$, $(x_1x_2,x_0x_2^2-x_1^3)$, if $A\cong A_{2,5}\oplus A_{0,1}$, $(x_2^2,x_1^2x_2-x_1^3)$, if $A\cong A_1^{sp}$, $(x_1x_2,x_2^3- x_1^3)$, if $A\cong A_2^{sp}$.
\qed
\enddemo

\subhead
6.2. The Scandinavian construction
\endsubhead
Let
$$
M:=\pmatrix
m_{0,0}&m_{0,1}&m_{0,2}\\
m_{1,0}&m_{1,1}&m_{1,2}\\
m_{2,0}&m_{2,1}&m_{2,2}
\endpmatrix
$$
be a matrix whose entries are linear forms $m_{i,j}(x_0,\dots,x_4)\in k[x_0,\dots,x_4]$, $0\le i,j\le2$. The ideal $F_2(M)\subseteq k[x_0,\dots,x_4]$ generated by the $2\times2$--minors of $M$ has depth at most $4$: if it is actually $4$ then its minimal free resolution is given by the Gulliksen--Neg\aa rd complex (see [G--N] or section 2.D of [B--V]: such a complex is called Scandinavian in [La]), thus $k[x_0,\dots,x_4]/F_2(M)$ is Gorenstein with Hilbert polynomial $p(t)=6$.

In particular the subscheme $X\subseteq\p4$ whose homogeneous ideal is $F_2(M)$ has dimension $0$, degree $6$ and it is arithmetically Gorenstein.

\definition{Definition 6.2.1}
The construction described above is called Scandinavian.
\enddefinition

Also the Scandinavian construction does not give all the points $X\in \Hilb_{6}^{aG,0}(\p{4})$. 
Indeed let $H\subseteq\p4$ be the linear subscheme whose ideal is generated by the entries of $M$. Clearly $H\subseteq X$, thus either $H=\emptyset$ or it consists of a single point $H=\{\ P\ \}$.

In the first case we can define an embedding
$\varphi\colon \p4\hookrightarrow\p8$ given by $x_{i,j}=m_{i,j}(x_0,\dots,x_4)$, $x_{i,j}$ being homogeneous coordinates in $\p8$ for $0\le i,j\le2$, thus $X=i^{-1}(S_{2,2})\cong S_{2,2}\cap\varphi(\p4)$, $S_{2,2}$ being the image of the standard Segre embedding $\p2\times\p2\hookrightarrow\p8$. Let $T_P(X)$ and $T_P(S_{2,2})$ be the tangent spaces at a point $P\in X$ and $P\in S_{2,2}$ respectively. Since $T_P(S_{2,2})\cap S_{2,2}$ is the union of two copies of $\p2\subseteq\p8$ intersecting only at $P$ and $\dim(S_{2,2}\cap\varphi(\p4))=0$, it turns out that $\dim(T_P(S_{2,2})\cap\varphi(\p4))\le2$, hence the same is true for $T_P(X)$.

In the second case, since the derivatives of the generators of the ideal of $X$ are in the ideal of $H$, then the tangent space at $P\in X$ is the whole of $\p4$, thus necessarily $X\in{\Cal G}_6$. E.g. let $i$ be any square root of $-1$ in $k$: then if
$$
M=\pmatrix
ix_4&x_1&x_2\\
-x_1&ix_4&x_3\\
-x_2&-x_3&ix_4
\endpmatrix
$$
then $X=G_6$.

\proclaim{Proposition 6.2.2}
$X\in \Hilb_{6}^{aG,0}(\p{4})$ can be obtained via the Scandinavian construction if and only if  $X\not\in O(A_{3,5}\oplus A_{0,1})\cup O(A_{3,6})$.
\endproclaim
\demo{Proof}
The statement is true if $A\cong A_{4,6}$ as we showed above.

In the other cases the idea is to consider a suitable embedding of $X$ in a complete intersection of two cubics in $\p2$ without common components, which intersect residually with respect to $X$ at three non--collinear point. Blowing out such three points we obtain $X$ as the intersection of the del Pezzo surface $S_6\subseteq\p6$ of degree $6$ with a linear space $H$ of dimension $4$. Since it is well--known the existence of a linear subspace $H'\subseteq \p8$ of dimension $2$ such that  $S_6=S_{2,2}\cap H'$, then $X=S_{2,2}\cap H$.

This is easy to do when $X\cong\spec(A)$ where $A$ is one of the following: $A_{1,2}^{\oplus e}\oplus A_{0,1}^{\oplus 6-2e}$ for $e=0,1,2,3$,  $A_{1,3}^{\oplus e}\oplus A_{0,1}^{\oplus 6-3e}$ for $e=1,2$,  $A_{1,3}\oplus A_{1,2}\oplus A_{0,1}$, $A_{1,4}\oplus A_{0,1}^{\oplus2}$, $A_{1,4}\oplus A_{1,2}$, $A_{1,5}\oplus A_{0,1}$, $A_{1,6}$. E.g. we now show how to obtain $X\cong\spec(A_{1,4}\oplus A_{1,2})$.

Fix a smooth cubic $F_1$ and choose a flex point $P_0\in F_1$ (with tangent line $T_{P_0}(F_1):=r_0$) and another general point $P_1\in F_1$ (with $T_{P_1}(F_1):=r_1$). Let $r$ be a general line through $P$ and $F_2:=r_0+r_1+r$. Then $F_1\cap F_2=X\cup Y$, where $X$ is a locally complete intersection (hence Gorenstein) scheme, whose associated cycle is $4P_0+2P_1$, wheras $Y$ does not intersect $X$ and it is supported on three pairwise distinct non--collinear points: hence $X\cong\spec(A_{n,4}\oplus A_{1,2})$ where $n$ is either $1$ or $2$. Since $T_{P_0}(X)\subseteq T_{P_0}(F_1)=r_0$, it follows $X\cong\spec(A_{1,4}\oplus A_{1,2})$.

Now we examine the remaining cases. For example consider two general conics $C_i$, $i=1,2$, intersecting at three points $P_0,P_1,P_2\in\p2$ such that $T_{P_0}(C_1)=T_{P_0}(C_2)$. Let $r_i$, $i=1,2$, be two general lines through $P_0$ and set $F_i:=C_i+r_i$, $i=1,2$. Then $F_1\cap F_2=X\cap Y$, where again $X$ is  Gorenstein with associated cycle $4P_0+2P_1$, $Y$ does not intersect $X$ and it is supported on three pairwise distinct non--collinear points. But since $T_{P_0}(X)=T_{P_0}(C_1)\cap T_{P_0}(C_2)$, it follows that $X\cong\spec(A_{2,4}\oplus A_{1,2})$ in this case. In a similar way one can handle the cases $A_{2,4}\oplus A_{0,1}^{\oplus2}$ and $A_{2,5}\oplus A_{0,1}$.

Let $F_1$ be an integral cubic with a cusp at $P_0$ and let $r_0$ be the unique tangent line to $F_1$ at $P_0$. Let $C$ be any integral conic through $P_0$ and tangent to $r_0$ and set $F_2:=C+r_0$. Then $F_1\cap F_2=X\cap Y$, where $X\cong\spec(A_1^{sp})$ and, for a general choice of $C$ the scheme $Y$ does not intersect $X$ and it is supported on three pairwise distinct non--collinear points.

Let $F_1$ be an integral cubic with a node at $P_0$ and let $r_0,r_1$ be the two tangent lines to $F_1$ at $P_0$. Let $C$ be any integral conic through $P_0$ and tangent to $r_0$ and set $F_2:=C+r_1$ (resp. $F_2:=C+r_0$). Then $F_1\cap F_2=X\cap Y$, where $X\cong\spec(A_2^{sp})$ (resp. $X\cong\spec(A_{2,6})$ and again, for a general choice of $C$ the scheme $Y$ does not intersect $X$ and it is supported on three pairwise distinct non--collinear points.
\qed
\enddemo

\remark{Remark 6.2.3}
In particular all the schemes that can be obtained via the Japanese construction can also be obtained via the Scandinavian construction.
\endremark
\medbreak

\subhead
6.3. The Anglo--American construction
\endsubhead
Let $M:=\pmatrix m_{i,j,h}\endpmatrix_{0\le i,j,h\le 1}$ be a cubic array whose entries are linear forms $m_{i,j,h}(x_0,\dots,x_4)\in k[x_0,\dots,x_4]$, $0\le i,j,h\le1$. Let $I\subseteq k[x_0,\dots,x_4]$ be the ideal generated by the determinants of all the faces of $M$ which are $12$ (we are taking into account also the diagonal faces). Such an ideal is generated by only $9$ of them, the other ones being linearly dependent, it has depth at most $4$ and, if it is $4$, then its minimal free resolution over $S:=k[x_0,\dots,x_4]$ looks like
$$
0\longrightarrow S(-6)\longrightarrow S(-4)^{\oplus9}\longrightarrow S(-3)^{\oplus16}\longrightarrow S(-2)^{\oplus9}\longrightarrow I\longrightarrow0.
$$
Thus $k[x_0,\dots,x_4]/I$ is Gorenstein with Hilbert polynomial $p(t)=6$. Such a resolution is obtained as a suitable mapping cone of a generalized Eagon--Northcott complex  (this is well--known: see section 2 of [Cs2] for a reference).

In particular the subscheme $X\subseteq\p4$ whose homogeneous ideal is $I$ has dimension $0$, degree $6$ and it is arithmetically Gorenstein.

\definition{Definition 6.3.1}
The construction described above is called Anglo--American.
\enddefinition

Let $H\subseteq\p4$ be the linear subscheme whose ideal is generated by the entries of $M$. Clearly $H\subseteq X$, thus either $H=\emptyset$ or it consists of a single point $H=\{\ P\ \}$.

In the first case we have an embedding
$\varphi\colon \p4\hookrightarrow\p7$ given by $x_{i,j,h}=m_{i,j,h}(x_0,\dots,x_4)$, $x_{i,j,h}$ being homogeneous coordinates in $\p7$ for $0\le i,j,h\le1$, thus $X=i^{-1}(S_{1,1,1})\cong S_{1,1,1}\cap\varphi(\p4)$, $S_{1,1,1}$ being the image of the standard Segre embedding $\p1\times\p1\times\p1\hookrightarrow\p7$. If $P\in X\subseteq S_{1,1,1}$, then $T_P(S_{1,1,1})\cap S_{1,1,1}$ is the union of three copies of $\p1\subseteq\p7$ intersecting only at $P$: since $\dim(S_{1,1,1}\cap\varphi(\p4))=0$, it turns out that $\dim(T_P(S_{1,1,1})\cap\varphi(\p4))\le2$, hence the same is true for $T_P(X)$.

In the second case, since the derivatives of the generators of the ideal of $X$ are in the ideal of $H$, then the tangent space at $P\in X$ is the whole of $\p4$, thus necessarily $X\in{\Cal G}_6$. E.g. if we choose
$$
\gather
m_{0,0,0}=-m_{1,1,1}=-x_3,\quad m_{0,0,1}=-m_{1,1,0}=x_1,\\
 m_{0,1,0}=m_{1,0,1}=ix_4,\quad m_{1,0,0}=-m_{0,1,1}=x_2
\endgather
$$
then $X=G_6$.

\proclaim{Proposition 6.3.2}
$X\in \Hilb_{6}^{aG,0}(\p{4})$ can be obtained via the Anglo--American construction if and only if it can be obtained via the Scandinavian construction.
\endproclaim
\demo{Proof}
The statement is true if $A\cong A_{4,6}$ as we showed above.

In the other cases it follows from the well--known fact that, for each general hyperplane $H\subseteq\p7$, the intersection $S_{1,1,1}\cap H$ is a del Pezzo surface $S_6\subseteq\p6$ of degree $6$.
\qed
\enddemo

\remark{Remark 6.3.3}
In the papers [Cs1] and [Cs2] some examples of particular constructions of covers of degree $6$ are described. The idea is to extend the above described Japanese, Scandinavian and Anglo--American construction to the case of schemes flat and finite of degree $6$ over an irreducible base $Y$. We will then speak respectively of Scandinavian, Anglo--American and Japanese covers.  

Due to numerical reasons it is proved that the classes of schemes $X$ that can be obtained as Scandinavian or Anglo--American covers of an irreducible $Y$ are not contained one into the other. It is not immediately clear that such covers are not Japanese.

Our study performed in Proposition 6.1.2, 6.2.2 and 6.3.2 and Example 4.6 of [Cs2] allows us to state that not all Scandinavian and Anglo--American covers can be obtained as Japanese covers since it is easy to construct examples of Scandinavian and Anglo--American covers $X\to Y$ with $X$ smooth and irreducible and having a fibre isomorphic to $\spec(A_{2,6})$, which is not allowed for Japanese covers.
\endremark
\medbreak

\subhead
6.4. The British construction
\endsubhead
We recall the definition of antisymmetric and extrasymmetric matrix format introduced firstly in Section 5.7 of [Rd1] (with another name). In a second time such a matrix format has been slightly generalized in [Rd2], [Rd3], [B--C--P1] and [B--C--P2].

\definition{Definition 6.4.1}
A $6\times6$ matrix $N$ with entries in a polynomial ring over a field $k$ is called antisymmetric and extrasymmetric if there are $3\times3$ matrices $A$ antisymmetric and $S$ symmetric and $q\in k$ such that 
$$
N=\pmatrix
A&S\\
-S&-qA
\endpmatrix.
$$
\enddefinition

Assume that the entries of $N$ are linear forms $n_{i,j}\in k[x_0,\dots,x_4]$, $0\le i,j\le 5$. Then for each general choice of the entries the saturated ideal of the locus in $\p4$ where the rank of $N$ is at most $2$ coincides with the ideal $Pf_4(N)$ generated by the pfaffians of order $4$ of $N$. A direct computation shows that such an ideal is generated by $9$ of the $15$ pfaffians. Moreover the syzygies among these pfaffians can be directly deduced by the matrix $N$ (see [Rd1], Section 5.9).

\definition{Definition 6.4.2}
The construction described above is called British.
\enddefinition

Let $X\subseteq\p4$ be the subscheme whose homogeneous ideal is $Pf_4(N)$. If $q\ne0$ a direct computation shows that $Pf_4(N)=F_2(M)$, $M$ being the matrix defined in section 6.2, thus $X\in \Hilb_{6}^{aG,0}(\p{4})$ for a general choice of $N$ and it can be also obtained via the Scandinavian construction.
If $q=0$ some easy computations show that if $\dim(X)=0$ then $X$ can be also obtained via the Japanese construction.

Taking into account Remark 6.2.3, we can summarize the above comments in the following 

\proclaim{Proposition 6.4.3}
$X\in \Hilb_{6}^{aG,0}(\p{4})$ can be obtained via the British construction if and only if it can be obtained via  the Scandinavian construction.
\qed
\endproclaim

\subhead
6.5. The Italian construction
\endsubhead
Consider an arithmetically Gorenstein subscheme $\widetilde{X}$ of degree $5$ contained in a hyperplane $h\subseteq\p4$, let $P\in\widetilde{X}$ . With a proper choice of the homogeneous coordinates $x_0,\dots,x_4$ in $\p4$ we can assume $h=\{\ x_4=0\ \}$, $P=[1,0,0,0,0]$.

The minimal free resolution of the ideal $I_P\subseteq S:=k[x_0,\dots,x_4]$ of $P$ is given by the Koszul complex.
In order to compute the minimal free resolution of the homogeneous ideal $I_{\widetilde{X}}\subseteq S:=k[x_0,\dots,x_4]$ of $\widetilde{X}$ we use the fact that $x_0,\dots,x_3$ are homogeneous coordinates in $h$, so a suitable mapping cone of the Buchsbaum--Eisenbud complex resolving the homogeneous ideal $I_{\widetilde{X}\vert h}\subseteq S:=k[x_0,\dots,x_3]$ yields the asked resolution (see [B--E]). 

Since $P\in X$ we can thus form a commutative diagram
$$
\matrix
0&\hskip-3truemm\longrightarrow&\hskip-3truemm S(-6)&\hskip-3truemm\longrightarrow&\hskip-3truemm\matrix S(-5)^{\phantom{\oplus5}}\\ \oplus^{\phantom{\oplus5}}\\ S(-4)^{\oplus5}\endmatrix&\hskip-3truemm\longrightarrow&\hskip-3truemm\matrix S(-3)^{\oplus5}\\ \oplus^{\phantom{\oplus5}} \\ S(-3)^{\oplus5}\endmatrix&\hskip-3truemm\longrightarrow&\hskip-3truemm\matrix S(-2)^{\oplus5}\\ \oplus^{\phantom{\oplus5}}\\ S(-1)^{\phantom{\oplus5}}\endmatrix&\hskip-3truemm\longrightarrow&\hskip-3truemm I_{\widetilde{X}}&\hskip-3truemm\longrightarrow&\hskip-3truemm0\\
&\hskip-3truemm &\hskip-3truemm\mapdown{\varphi}&\hskip-3truemm &\hskip-3truemm\mapdown{\phantom{\varphi}}&\hskip-3truemm &\hskip-3truemm\mapdown{\phantom{\varphi}}&\hskip-3truemm &\hskip-3truemm\mapdown{\phantom{\varphi}}&\hskip-3truemm &\hskip-3truemm\mapdown{\phantom{\varphi}}&\hskip-3truemm &\hskip-3truemm\\
0&\hskip-3truemm\longrightarrow&\hskip-3truemm S(-4)&\longrightarrow&\hskip-3truemm S(-3)^{\oplus4}&\hskip-3truemm\longrightarrow&\hskip-3truemm S(-2)^{\oplus6}&\hskip-3truemm \longrightarrow&\hskip-3truemmS(-1)^{\oplus4}&\hskip-3truemm\longrightarrow&\hskip-3truemm I_P&\hskip-3truemm\longrightarrow&\hskip-3truemm 0
\endmatrix
$$
where $\varphi$ is represented by a single quadratic form $f\in k[x_0,\dots,x_3]$. Due to Theorem 3.2 of [B--D--N--S] the ring $S_X:=k[x_0,\dots,x_4]/I_X$ where
$$
I_X:=I_{\widetilde{X}\vert h}+(x_1x_4,x_2x_4,x_3x_4, f+x_4g)\subseteq k[x_0,\dots,x_4]
$$
is Gorenstein of codimension $4$ for each linear form $g\in k[x_0,\dots,x_4]\setminus k[x_0,\dots,x_3]$ (since in this case the sequence $x_1,x_2,x_3, f+x_4g$ is regular), thus the corresponding scheme $X\subseteq\p4$ is arithmetically Gorenstein. Moreover $\deg(X)=6$ due to Corollary 3.13 of [B--D--N--S].

Notice that the line $r:=\{\ x_1=x_2=x_3=0\ \}$ intersects $X$ along the scheme $X\cap r$ whose homogeneous ideal is $I_{X\cap r}=(x_1,x_2,x_3, f+x_4g)$, which has degree $2$. Let $Z$ be the residual scheme inside $\widetilde{X}$ to $P$: $h$ is  the unique hyperplane containing it (otherwise there would exist a hyperplane section of $X$ of degree $5$, an absurd by [Sch], Lemma (4.2)): since $I_X\subseteq I_{Z}$ and $I_X:I_{Z}=I_{X\cap r}$ by Proposition 3.15 of  [B--D--N--S], thus ${Z}$ and $X\cap r$ are $G$--linked via $X$, hence $\deg(Z)=4$. Since $X$ is arithmetically Gorenstein in $\p4$, then $\deg(X\cap h)\le4$, thus $Z=X\cap h$ due to the obvious inclusion $Z\subseteq X\cap h$.

\definition{Definition 6.5.1}
The construction described above is called Italian.
\enddefinition

We have proved the first part of the following characterization.

\proclaim{Lemma 6.5.2}
$X\in  \Hilb^{6}_{aG,0}(\p{4})$ can be obtained via the Italian construction if and only if there are a line $r$ and a hyperplane $h$ of the projective space $\p4$ satisfying the following three conditions:
\item{i)} $r\not\subseteq h$;
\item{ii)} $\deg(X\cap r)=2$;
\item{iii)} the residual scheme to $X\cap r$ inside $X$ is $X\cap h$.
\endproclaim
\demo{Proof}
We have only to prove the converse, so we assume the existence of $r,h\subseteq\p4$ and we will construct $P$ and $\widetilde{X}$ giving $X$ via the Italian construction. 

With a proper choice of the homogeneous coordinates $x_0,\dots,x_4$ in $\p4$, we can assume that $r:=\{\ x_1=x_2=x_3=0\ \}$ and $h:=\{\ x_4=0\ \}$ (condition i)). Thus $I_{X\cap r}=(x_1,x_2,x_3,q)$ (condition ii)) and $I_X=(q_1,\dots,q_9)$ for suitable quadratic forms $q,q_1,\dots,q_9\in k[x_0,\dots,x_4]$.

Since $I_X\subseteq I_{X\cap r}$ and $I_X\not\subseteq (x_1,x_2,x_3)$ (because $r\not\subseteq X$) it follows that $q$ can be chosen as a minimal generator of $I_X$, hence we can write $I_X=(q_1,\dots,q_8,q)$ with $q_i$ satisfying the conditions of Theorem 3.2 of  [B--D--N--S] and such $q_1,\dots,q_8\in k[x_1,x_2,x_3]$, thus we can apply the construction described there with $I=I_X$ and $J=(x_1,x_2,x_3,q)$ obtaining an arithmetically Gorenstein scheme $\widetilde{X}\subseteq\p4$ of degree 5 (again apply Corollary 3.13) whose homogeneous ideal is $I_{\widetilde{X}}:=(q_1,\dots,q_8,x_1q,x_2q,x_3q,\ell)\subseteq k[x_0,\dots,x_4]$ where $\ell$ is a linear form representing the map $\alpha_c$ in diagram (3) of [B--D--N--S]. It follows that $I_{\widetilde{X}}=(g_1,\dots,g_5,\ell)$.

By  Lemma 2.9 of [B--D--N--S], $\{\ \ell=0\ \}$ contains the residual scheme in $X$ with respect to $X\cap r$. Since $h$ is the unique hyperplane containing $X\cap h\subseteq\widetilde{X}$ (due again to [Sch], Lemma (4.2)), then condition iii) yields $\ell=x_4$, thus $I_{\widetilde{X}}\subseteq I_P:=(x_1,x_2,x_3,x_4)$. 
Now Proposition 3.17 of [B--D--N--S] guarantees that the Italian construction from the data $P\in\widetilde{X}$ yields the ideal $I_X$.
\qed
\enddemo

As for the previous constructions we now list all the possible Gorenstein algebras $A$ such that $X:=\spec(A)$ can be obtained via the Italian one.

\proclaim{Proposition 6.5.3}
$X\in \Hilb_{6}^{aG,0}(\p{4})$ can be obtained via the Italian construction  if and only if  $X\not\in O(A_{1,6})\cup O(A_1^{sp})\cup O(A_2^{sp})$.
\endproclaim
\demo{Proof}
First we examine the case when $X\in  \Hilb_{6}^{aG,0}(\p{4})$ is irreducible. In this case $I_X$ is given by Lemma 5.2. 

If $A\cong A_{2,6}$ then take $I_{\widetilde{X}}:=(x_1x_2-x_0x_3,x_2x_3,x_1^2-x_0x_2,x_2^2-x_4^2,x_3^2,x_4)$: thus $I_{\widetilde{X}}:I_P=I_{\widetilde{X}}+(x_1x_3)$ hence $f=x_1x_3$ (Lemma 2.9 of [B--D--N--S]). We conclude that the choice $g=-x_4$ yields that such a scheme can be obtained via the Italian construction.

If $A\cong A_{3,6}$ then take $I_{\widetilde{X}}:=(x_1x_3,x_2x_3,x_1^2-x_0x_2,x_2^2,x_3^2-x_4^2,x_4)$: thus $I_{\widetilde{X}}:I_P=I_{\widetilde{X}}+(x_1x_2)$ hence $f=x_1x_2$ and again $g=-x_4$ yields that such a scheme can be obtained via the Italian construction.

If $A\cong A_{4,6}$ then $I_{\widetilde{X}}:=(x_1x_2,x_1x_3,x_2x_3,x_2^2-x_1^2,x_3^2-x_1^2,x_4)$ and $I_{\widetilde{X}}:I_P=I_{\widetilde{X}}+(x_1x_2)$ hence $f=x_1x_2$ and again $g=-x_4$ yields that also such a scheme can be obtained via the Italian construction.

We now show there are no pairs $(\alpha,r)$ as in Lemma 6.5.2 when $A\cong A_{1,6},A_1^{sp},A_2^{sp}$. Let us consider the general line $r_{[a_1,\dots,a_4]}=\{\ [s,a_1t,\dots,a_4t]\ \}$ through $P$.
In any case if $a_2a_3a_4\ne0$ it is easy to check that $\deg(X\cap r_{[a_1,\dots,a_4]})=1$.

If $A\cong A_{1,6}$ the same is false only for $r_{[1,0,0,0]}$. In this case the ideal of $Y:=X\cap r_{[1,0,0,0]}$ is $I_Y=(x_2,x_3,x_4,x_1^2-x_0x_2)$, thus a direct computation shows that $x_4\in I_X:I_Y$. Since $x_4\in I_Y$ it follows that such an $X$ cannot be obtained via the Italian construction.

If $A\cong A_{1}^{sp}, A_2^{sp}$ the lines for which $\deg(X\cap r_{[a_1,\dots,a_4]})=2$ are $r_{[a_1,0,a_3,0]}$, then the ideal of $Y:=X\cap r_{[a_1,0,a_3,0]}$ is $I_Y=(a_3x_1-a_1x_3,x_2,x_4,x_3^2)$. In both the cases $a_3x_2+a_1x_4\in(I_X:I_Y)\cap I_Y$, hence we can prove as above that also in this case $X$ cannot be obtained via the Italian construction.

We now consider the reducible cases. The case $A\cong A_{0,1}^{\oplus6}$ has been already described in Example 3.18 of [B--D--N--S] where it is proved it can be obtained via the Italian construction.

If $A\cong A_{1,2}^{\oplus e}\oplus A_{0,1}^{\oplus 6-2e}, A_{1,3}\oplus A_{1,2}\oplus A_{0,1},A_{1,4}\oplus A_{1,2},A_{2,4}\oplus A_{1,2}$ ($e=1,2,3$) we take as $r$ any line generated by any subscheme $X\cong\spec(A)\subseteq\p4$ corresponding to a summand $A_{1,2}$. If $A\cong A_{1,3}\oplus A_{0,1}^{\oplus 3},A_{1,4}\oplus A_{0,1}^{\oplus2},A_{2,4}\oplus A_{0,1}^{\oplus2}$ we take any line $r$ through two simple points of the embedded scheme $X\cong\spec(A)\subseteq\p4$. In all the above cases the hyperplane $h$ containing the residual scheme in $X$ with respect to $X':=X\cap r$ cannot contain $r$, otherwise $\deg(X\cap h)\ge5$, an absurd by [Sch], Lemma (4.2).

Consider the case $A\cong A_{1,3}^{\oplus 2}$: in this case we have a decomposition in disjoint subschemes $X:=X_1\cup X_2\subseteq\p4$ where $X_i\cong\spec(A_{1,3})$. Let $\alpha_1\subseteq\p4$ the linear space generated by $X_i$: due to [Sch], Lemma (4.2) it is a plane. It follows the existence of a line $r\subseteq\p4$ such that $\deg(X_2\cap r)=2$. We have $r\cap \alpha_1=\emptyset$ otherwise the hyperplane generated by $\alpha_1$ and $r$ would intersect $X$ in a subscheme of degree $5$, an absurd. Let $h\subseteq\p4$ be the hyperplane generated by $\alpha_1$ and the support of $X_2$ which contains obviously the residual scheme with respect to $X\cap r$: by the above discussion it follows that $r\not\subseteq h$. 

Let $A\cong A_{3,5}\oplus A_{0,1}$. To this purpose we first check that the subscheme $Y\subseteq\p4$ whose homogeneous ideal in $S:= k[x_0,x_1,x_2,x_3,x_4]$ is
$$
I=(x_0x_1,x_0x_2,x_0x_3,x_1x_2,x_1x_3,x_2x_3,x_1^2-x_0x_4,x_2^2-x_0x_4,x_3^2-x_0x_4),
$$
satisfies $Y\in O(A)$.

It is easy to check that $Y=Y'\cup Y''$ where the supports of $Y'$ and $Y''$ are $P':=[1,0,0,0,0]$ and $P'':=[0,0,0,0,1]$ respectively. Since 
$$
\gather
I_{Y\cap\{\ x_4\ne0\ \}}=I_{Y'\cap\{\ x_4\ne0\ \}}=(x_1x_2,x_1x_3,x_2x_3,x_1^2-x_2^2, x_1^2-x_3^2,x_3^2-x_0),\\
I_{Y\cap\{\ x_0\ne0\ \}}=I_{Y''\cap\{\ x_0\ne0\ \}}=(x_1,x_2,x_3,x_4)
\endgather
$$
one obtains the isomorphism $Y\cong\spec(A)$. Moreover one also checks in this way that $\deg(Y\cap \alpha)\le4$ for each hyperplane $\alpha\subseteq\p4$, thus the embedding of $Y$ is arithmetically Gorenstein by Lemma (4.2) of [Sch].

Thus we can assume that the ideal of $X$ is actually $I_X=I$. Consider the line $r:=\{\ x_1=x_2=x_3=0\ \}$ and the hyperplane $h:=\{\ x_0=0\ \}$. Thus $r\not\subseteq h$, $\deg(X\cap r)=2$ and $I_X:I_{X\cap r}\supseteq I_X+(x_0)$, thus the residual scheme $Z$ in $X$ with respect to $X\cap r$ is contained in the scheme $X\cap h$ whose homogeneous ideal is $I_X+(x_0)$. Since $I_X+(x_0)=(x_1,x_2,x_3)^2+(x_0)$, trivially $\deg(Z)=4$, whence it follows that $Z=X\cap h$. We conclude that such an $X$ can be obtained via the Italian construction by Lemma 6.5.2.

If  either $A\cong A_{2,5}\oplus A_{0,1}$ or $A\cong A_{1,5}\oplus A_{0,1}$, then we repeat the above argument with 
$$
\gather
I:=(x_1x_2,x_1x_3-x_0x_4,x_1x_4,x_2x_3,x_2x_4,x_3x_4,x_1^2-x_0x_3,x_2^2-x_0x_4,x_3^2),\\
I:=(x_1x_2-x_0x_3,x_1x_3-x_0x_4,x_1x_4,x_2x_3,x_2x_4,x_3x_4,x_1^2-x_0x_2,x_2^2-x_0x_4,x_3^2)
\endgather
$$
respectively.
\qed
\enddemo

\medbreak

\subhead
6.6. The Anglo--Hellenic construction
\endsubhead
The following construction is based on the Tom unprojection described in Section 5 of [Pa] (about unprojections see also the fundamental papers [P--R] and [K--M]).

Let 
$$
A:=\pmatrix
0&a_{0,1}&a_{0,2}&a_{0,3}&a_{0,4}\\
-a_{0,1}&0&a_{1,2}&a_{1,3}&a_{1,4}\\
-a_{0,2}&-a_{1,2}&0&a_{2,3}&a_{2,4}\\
-a_{0,3}&-a_{1,3}&-a_{2,3}&0&a_{3,4}\\
-a_{0,4}&-a_{1,4}&-a_{2,4}&-a_{3,4}&0
\endpmatrix
$$
be a skew--symmetric matrix with linear entries $a_{i,j}\in k[x_0,\dots,x_4]$, $0\le i<j\le4$. We make two assumptions from now on. First we assume that the locus where the pfaffians of order $4$ of $A$ vanish is a curve $\widetilde{X}\subseteq\p4$. Second we assume that $a_{i,j}\in k[x_1,\dots,x_4]$, $1\le i<j\le4$.

Each such matrix $A$ gives rise, via unprojection, to a new scheme $\widetilde{X}_{unp}\subseteq\p5$ which is again an arithmetically Gorenstein curve (see Theorem 1.5 of [P--R]), thus each hyperplane $h\subseteq\p5$ not containing any of its component cuts $\widetilde{X}_{unp}$ along a $0$--dimensional scheme $X:=\widetilde{X}_{unp}\cap h\subseteq h\cong\p4$ which is arithmetically Gorenstein too.

For each $i=0,\dots,4$ we denote by $p_i$ the pfaffian of the submatrix $A_i$ obtained from $A$ by deleting the $i^{\roman th}$ row and column. Since $a_{i,j}\in (x_1,\dots,x_4)$, $1\le i<j\le4$, then we can find a $4\times4$ matrix $Q$ whose entries are linear combinations of $a_{0,j}$, $j=1,\dots,4$ such that $(p_1,\dots,p_4)=(x_1,\dots,x_4)Q$. Then  the cofactors of the first column of $Q$ are multiple of $a_{0,1}$: let $g_1,\dots,g_4$ be such cofactors divided by $a_{0,1}$.

Then
$$
I_{\widetilde{X}_{unp}}=I_{\widetilde{X}}+(x_1s-g_1,x_2s-g_2,x_3s-g_3,x_4s-g_4)\subseteq k[x_0,\dots,x_4,s].
$$
is the ideal of $\widetilde{X}_{unp}$. Since $x_jg_i-x_ig_j\in I_{\widetilde{X}}$ by Lemma 5.4 of [Pa], a Groebner basis of $I_{\widetilde{X}_{unp}}$ is given by a Groebner basis of $I_{\widetilde{X}}$ plus the four polynomials  $x_is-g_i$. It follows that $I_{\widetilde{X}_{unp}}\cap k[x_0,\dots,x_4]=I_{\widetilde{X}}$, i.e. the projection of $\widetilde{X}_{unp}$ form $Q:=[0,0,0,0,0,1]\in \widetilde{X}_{unp}$ is exactly $\widetilde{X}$. Since $\deg(\widetilde{X})=5$ and $Q$ is simple on $\widetilde{X}_{unp}$ then $\deg(\widetilde{X}_{unp})=6$, then $X=\widetilde{X}_{unp}\cap h\in  \Hilb_{6}^{aG,0}(\p{4})$. 

\definition{Definition 6.6.2}
The construction described above is called Anglo--Hellenic.
\enddefinition

The following Proposition shows that the Anglo--Hellenic construction is actually a Babylonian construction.

\proclaim{Proposition 6.6.3}
Every $X\in \Hilb_{6}^{aG,0}(\p{4})$ can be obtained via the Anglo--Hellenic construction.
\endproclaim
\demo{Proof}
Assume that $X$ can be obtained as the degeneracy locus of the matrix
$$
M:=\pmatrix
m_{0,0}&m_{0,1}&m_{0,2}\\
m_{1,0}&m_{1,1}&m_{1,2}\\
m_{2,0}&m_{2,1}&m_{2,2}
\endpmatrix.
$$
It is then proved in section 5.5 of [Pa] that it suffices to unproject the scheme $\widetilde{X}$ associated to the skew--symmetric matrix 
$$
A=\pmatrix
0&m_{1,0}&m_{2,0}&m_{0,1}&m_{0,2}\\
-m_{1,0}&0&0&m_{1,1}&m_{1,2}\\
-m_{2,0}&0&0&m_{2,1}&m_{2,2}\\
-m_{0,1}&-m_{1,1}&-m_{2,1}&0&0\\
-m_{0,2}&-m_{1,2}&-m_{2,2}&0&0
\endpmatrix
$$
and to intersect with $h:=\{\ s-m_{0,0}=0\ \}$, i.e. each $X\in \Hilb_{6}^{aG,0}(\p{4})$ which can be obtained via the Scandinavian construction, can also be obtained via the Anglo--Hellenic one. 

Then we have only to check the statement for $X\in O(A_{3,5}\oplus A_{0,1})\cup O(A_{3,6})$. As we showed in the proof of Proposition 6.5.3, if $X\in O(A_{3,5}\oplus A_{0,1})$, then up to projectivities
$$
I_X=(x_0x_1,x_0x_2,x_0x_3,x_1x_2,x_1x_3,x_2x_3,x_1^2-x_0x_4,x_2^2-x_0x_4,x_3^2-x_0x_4).
$$
It is easy to check that $I_X$ can be obtained via the above construction with 
$$
A=\pmatrix
0&-x_3&-x_1&x_1&x_2\\
x_3&0&0&x_2&x_1\\
x_1&0&0&0&x_3\\
-x_1&-x_2&0&0&0\\
-x_2&-x_1&x_3&0&0
\endpmatrix
$$
and $s=x_0$.

By Lemma 5.2 iii), if $X\in O(A_{3,6})$, then up to projectivities
$$
I_X=(x_1x_2-x_4^2,x_1x_3,x_1x_4,x_2x_3,x_2x_4,x_3x_4,x_1^2-x_0x_2,x_2^2,x_3^2-x_4^2).
$$
Again $I_X$ can be obtained by choosing
$$
A=\pmatrix
0&x_2&x_1&x_3&x_0\\
-x_2&0&0&0&x_3\\
-x_1&0&0&-x_2&0\\
-x_3&0&x_2&0&x_1\\
-x_0&-x_3&0&-x_1&0
\endpmatrix
$$
and $s=x_4$.
\qed
\enddemo

\head
7. From $\Hilb_d^{aG,0}(\p{d-2})$ to $\Hilb_{d+1}^{aG,0}(\p{d-1})$
\endhead

We already know that the singular locus $\Sing(\Hilb^{aG,0}_{d}(\p{d-2}))$ of $\Hilb^{aG,0}_{d}(\p{d-2})$ is non--empty since it contains the locus ${\Cal G}_d$ for $d\ge6$ (see Theorem 3.5). Moreover we also checked that equality holds if $d=6$. It is then natural to ask if there exists other obstructed schemes in $\Hilb^{aG,0}_{d}(\p{d-2})$ as $d$ grows up. The aim of this section is to answer affirmatively this question.

\proclaim{Lemma 7.1}
Let $X\in\Hilb^{aG,0}_{d}(\p{d-2})$ be a scheme having a reduced component $P$. Then the restriction of the projection with center $P$ onto a hyperplane $H\cong\p{d-3}$ induces an isomorphism between $X\setminus P$ and its image $X_P\subseteq H$. Moreover $X_P\in\Hilb^{aG,0}_{d-1}(H)$.
\endproclaim
\demo{Proof}
Let $\pi\colon \p{d-2}\dashrightarrow H$. In order to prove that $\pi_{\vert X\setminus P}$ is an isomorphism onto its image it suffices to check that the secant and the tangent varieties of $X\setminus P$ do not contain $P$. Let $x',x''\in X$ (possibly coinciding) points. If $x',x'',P$ were collinear then the line $r$ through $x'$ and $x''$ would be contained in each quadric through $X$. Since $X$, being arithmetically Gorenstein and non--degenerate, is intersection of quadrics due to Sequence (2.3.1), this yields an absurd.

Since $X$ is non--degenerate then the same is true for $X_P$. Now let $\alpha\subseteq H$ be a hyperplane such that $\deg(X_P\cap \alpha)\ge d-2$: then the hyperplane $\beta\subseteq\p{d-2}$ generated by $\alpha$ and $P$ would be such that $\deg(\beta\cap X)\ge d-1$, an absurd. This last remark concludes the proof due to [Sch], Lemma (4.2).
\qed
\enddemo

Now let $X\in\Hilb^{aG,0}_{d}(\p{d-2})$ containing a reduced component $P$. In order to compute $h^0\big(X,{\Cal N}_X\big)$ consider a quadric $Q\subseteq\p{d-2}$ containing $X$. Since  $\Hom_{\p{d-2}}\big({\Cal F},\Ofa{X}\big)\cong\Hom_{X}\big({\Cal F}\otimes\Ofa X,\Ofa{X}\big)$ for each sheaf $\Cal F$ of $\Ofa X$--modules, then we have the exact sequence
$$
0\longrightarrow{\Cal N}_{X\vert Q}\longrightarrow{\Cal N}_{X}\longrightarrow{\Cal N}_{Q}\otimes\Ofa{X}\longrightarrow\Ext_{X}^1\big(\Im_{X\vert Q}/\Im_{X\vert Q}^2,\Ofa{X}\big),
$$
where $\Im_{Y\vert Q}\subseteq\Ofa Q$ denotes the sheaf of ideals of the subscheme $Y\subseteq Q$ and ${\Cal N}_{Y\vert Q}:=\Hom_{Q}\big(\Im_{Y\vert Q}/\Im_{Y\vert Q}^2,\Ofa{Y}\big)$. Since $X$ is locally Gorenstein then $\Ext_{X}^1\big(\Im_{X\vert Q}/\Im_{X\vert Q}^2,\Ofa{X}\big)=0$. Since all the sheaves are supported on the affine scheme $X$ and ${\Cal N}_{Q}\cong\Ofa Q(2)$, taking cohomologies we obtain
$$
h^0\big(X,{\Cal N}_{X}\big)= h^0\big(X,{\Cal N}_{X\vert Q}\big)+h^0\big(X,\Ofa X(2)\big)=h^0\big(X,{\Cal N}_{X\vert Q}\big)+d.
$$
We can assume that $Q$ is smooth thanks to the following

\proclaim{Lemma 7.2}
Let $X\in\Hilb^{aG,0}_{d}(\p{d-2})$, $d\ge4$. Then the general quadric through $X$ is smooth if and only if $X\not\in{\Cal G}_d$.
\endproclaim
\demo{Proof}
The homogeneous ideal $I_X\subseteq k[x_0,\dots,x_{d-2}]$ is generated by quadrics (see Sequence (2.3.1)). In particular the general quadric through $X$ is smooth outside $X$.

Since the tangent space $T_x(X)$ at the point $x\in X$ is exactly $\bigcap T_x(Q)$ where $Q$ runs in the set of quadrics through $X$, it follows that there are no smooth quadrics through $X$ if and only if there is $x\in X$ such that $T_x(X)=\p{d-2}$, i.e. $X\in{\Cal G}_d$.
\qed
\enddemo

We have a natural decomposition $h^0\big(X,{\Cal N}_{X\vert Q}\big)=h^0\big(X\setminus P,{\Cal N}_{X\setminus P\vert Q}\big)+h^0\big(P,{\Cal N}_{P\vert Q}\big)$. Since $Q$ is smooth then $h^0\big(P,{\Cal N}_{P\vert Q}\big)=d-3$. For a general choice of $Q$ we can assume that $X\setminus P\subseteq Q\setminus T_P(Q)$. The projection $\pi$ from $P$ onto $H$ induces an isomorphism $\pi'$ from $Q\setminus T_P(Q)$ onto $H\setminus T_P(Q)$, thus $h^0\big(X\setminus P,{\Cal N}_{X\setminus P\vert Q}\big)=h^0\big(X_P,{\Cal N}_{X_P}\big)$. We thus conclude
$$
h^0\big(X,{\Cal N}_{X}\big)=h^0\big(X_P,{\Cal N}_{X_P}\big)+2d-3.\tag7.3
$$

We are now ready to state the main result of this section. To this purpose we denote by ${O_{h,d}}$ the closure of $O(A_{h,h+2}\oplus A_{0,1}^{\oplus d-h-2})$ inside $\Hilb^{aG}_{d}(\p{d-2})$.

\proclaim{Theorem 7.4}
If $d\ge6$ then ${O_{4,d}}\subseteq\Sing(\Hilb^{aG}_{d}(\p{d-2}))$.
\endproclaim
\demo{Proof}
First we prove $O_{4,d}\subseteq\Hilb^{gen}_{d}(\p{d-2})$. Since, by definition, the last scheme is closed in $\Hilb^{aG}_{d}(\p{d-2})$, it suffices to check $O(A_{4,6}\oplus A_{0,1}^{\oplus d-6})\subseteq\Hilb^{gen}_{d}(\p{d-2})$. Since $\Hilb^{gen}_{d}(\p{d-2})$ is the closure of $O(A_{0,1}^{\oplus d})$, it suffices to prove that $O(A_{4,6}\oplus A_{0,1}^{\oplus d-6})\cap \Hilb^{gen}_{d}(\p{d-2})\ne\emptyset$.

We prove the assertion  by induction on $d\ge 6$, the first step being trivial due to Corollary 3.4. Assume we have proved the assertion for some $d-1\ge 6$ and consider $X\in O(A_{4,6}\oplus A_{0,1}^{\oplus d-6})$: since the projection $\pi$ from each simple point $P\in X$ onto $H$ induces, for each quadric $Q$ smooth along $X$, an isomorphism $\pi'\colon Q\setminus T_P(Q)$ onto $H\setminus T_P(Q)$ then $X_P\in O(A_{4,6}\oplus A_{0,1}^{\oplus d-7})$. Since $O(A_{4,6}\oplus A_{0,1}^{\oplus d-7})\subseteq\Hilb^{gen}_{d-1}(\p{d-3})$ which is the closure of $O(A_{0,1}^{\oplus d-1})$ by induction hypothesis we can find a flat family ${\Cal X}'\to B'$ with special fibre isomorphic to $X_P$ at $b_0\in B'$ and reduced general fibre. Over a suitable open neighborhood $B\subseteq B'$ of $b_0$, the inverse images of the fibres of ${\Cal X}'$ via $\pi'$ give rise to a flat family ${\Cal X}\to B$ such that ${\Cal X}\cup P\to B$ is again flat and it has special fibre isomorphic to $X$ over $b_0$ and reduced general fibre. In particular $X\in O(A_{4,6}\oplus A_{0,1}^{\oplus d-6})\cap \Hilb^{gen}_{d}(\p{d-2})$.

In particular we have proved that $X\in O(A_{4,6}\oplus A_{0,1}^{\oplus d-6})$ are all contained in an irreducible component of dimension $\dim(\Hilb^{gen}_{d}(\p{d-2}))=d(d-2)$. Now we prove that each $X\in O(A_{4,6}\oplus A_{0,1}^{\oplus d-6})$ is obstructed again by induction on $d\ge 6$, the first step being Theorem 3.5. Assume the statement true for some $d-1\ge 6$. If $X\in O(A_{4,6}\oplus A_{0,1}^{\oplus d-6})$ projecting from a simple point $P\in X$ we obtain $X_P\in O(A_{4,6}\oplus A_{0,1}^{\oplus d-7})$, which is obstructed by induction hypothesis. Thanks to Formula (7.3) we thus obtain $h^0\big(X,{\Cal N}_{X}\big)>(d-1)(d-3)+2d-3=d(d-2)$, whence $X$ is obstructed too.
\qed
\enddemo

We know that ${\Cal G}_d=O_{d-2,d}$ hence $\Sing(\Hilb^{aG}_{d}(\p{d-2}))$ contains the two closed subschemes $O_{d-2,d}, O_{4,d}$. It is then natural to deal with the intermediate cases ${O_{h,d}}$ when $h>4$. 

\proclaim{Proposition 7.5}
For each $h\ge0$ there exists an inclusion $O_{h+1,d}\subseteq O_{h,d}$. In particular ${O_{h,d}}\subseteq\Sing(\Hilb^{aG}_{d}(\p{d-2}))$ when $h\ge4$.
\endproclaim
\demo{Proof}
Since the $O_{h+1,d}$ is the closure of the orbit $O(A_{h+1,h+3}\oplus A_{0,1}^{\oplus d-h-3})$ it suffices to check $O(A_{h+1,h+3}\oplus A_{0,1}^{\oplus d-h-3})\cap O_{h,d}\ne \emptyset$. To this purpose we have to find a family with special fibre in $O(A_{h+1,h+3}\oplus A_{0,1}^{\oplus d-h-3})$ and general fibre in $O(A_{h,h+2}\oplus A_{0,1}^{\oplus d-h-2})$. If we construct a family with this property in the case $d=h+2$, then adding simple points we obtain the statement for each $d$.

Take $X\in  O(A_{d-3,d}\oplus A_{0,1})$: as we showed in Proposition 3.3 there are homogeneous coordinates $x_0,\dots,x_{d-2}$ such that the homogeneous ideal $I_X\subseteq k[x_0,\dots,x_{d-2}]$ of $X$ is
$$
I_X=(x_ix_j-\delta_{i,j}x_1^2+x_0\ell_{i,j})_{1\le i\le j\le {d-2},\ (i,j)\ne(1,1)}
$$
where $\ell_{i,j}\in k[x_1,\dots,x_{d-2}]$ are linear forms for $i,j=1,\dots,{d-2}$. Consider the family
$$
{\Cal A}:=k[b,x_0,\dots,x_{d-2}]/(x_ix_j-\delta_{i,j}x_1^2+bx_0\ell_{i,j})_{1\le i\le j\le {d-2},\ (i,j)\ne(1,1)}\to\a1.
$$
Notice that if $b\ne0$ the transformation $x_0\mapsto x_0/b$ induces an isomorphism $\spec({\Cal A}_b)\cong X\in O(A_{d-3,d}\oplus A_{0,1})$, while $\spec({\Cal A}_0)\in{\Cal G}_d$. 
\qed
\enddemo

\remark{Remark 7.6}
We have proved above the existence of a chain$$
{\Cal G}_d=O_{d-2,d}\subseteq O_{d-3,d}\subseteq \dots\subseteq  O_{5,d}\subseteq O_{4,d}\subseteq\Sing(\Hilb^{aG}_{d}(\p{d-2})),
$$
thus the following questions arise naturally.
\definition{Question 7.6.1}
Is $O_{4,d}=\Sing(\Hilb^{aG}_{d}(\p{d-2}))$? If not, does there exist a geometric description of $\Sing(\Hilb^{aG}_{d}(\p{d-2}))\setminus O_{4,d}$?
\enddefinition

It is clear that each $X\in O_{h,d}$ contains a point $P$ such that $\dim(T_P(X))\ge h$. Moreover if $X\in O(A_{h,h+2}\oplus A_{0,1}^{\oplus d-h+2})$ then it contains a component isomorphic to $\spec(A_{h,h+2})$.
\definition{Question 7.6.2}
Does $O_{h,d}$ coincide with the set of $X\in\Hilb^{aG}_{d}(\p{d-2})$ having tangent space of dimension at least $h$ at some point? Does $O_{h,d}$ coincide with the set of $X\in\Hilb^{aG}_{d}(\p{d-2})$ containing $\spec(A_{h,h+2})$ as subscheme?
\enddefinition
\endremark
\medbreak

\enddocument